\documentclass[11pt]{article}

\usepackage{amsmath}  
\usepackage{amsthm} 
\usepackage{amssymb} 
\usepackage[backref=page]{hyperref}  
\usepackage{enumitem} 
\usepackage{color}

\newtheorem{thm}{Theorem}
\newtheorem{inspr}[thm]{}

\newenvironment{definitie}{\begin{itemize}\item[ ]\hspace{-26pt}\bf Definition \rm }{\end{itemize}}
\newenvironment{notatie}{\begin{itemize}\item[ ]\hspace{-26pt}\bf Notation \rm }{\end{itemize}}
\newenvironment{voorbeeld}{\begin{itemize}\item[ ]\hspace{-26pt}\bf Example \rm }{\end{itemize}}
\newenvironment{stelling}{\begin{itemize}\item[ ]\hspace{-26pt}\bf Theorem \rm }{\end{itemize}}
\newenvironment{propositie}{\begin{itemize}\item[ ]\hspace{-26pt}\bf Proposition \rm }{\end{itemize}}
\newenvironment{lemma}{\begin{itemize}\item[ ]\hspace{-26pt}\bf Lemma \rm }{\end{itemize}}
\newenvironment{opmerking}{\begin{itemize}\item[ ]\hspace{-26pt}\bf Remark \rm }{\end{itemize}}
\newenvironment{voorwaarde}{\begin{itemize}\item[ ]\hspace{-26pt}\bf Condition \rm }{\end{itemize}}
\newenvironment{probleem}{\begin{itemize}\item[ ]\hspace{-26pt}\bf Problem \rm }{\end{itemize}}
\newenvironment{gevolg}{\begin{itemize}\item[ ]\hspace{-26pt}\bf Corollary \rm }{\end{itemize}}
\newenvironment{niets}{\begin{itemize}\item[ ]\hspace{-26pt}\bf   \rm }{\end{itemize}}

\newcommand{\defin}{\begin{inspr}\begin{definitie}}  
\newcommand{\edefin}{\end{definitie}\end{inspr}}

\newcommand{\notat}{\begin{inspr}\begin{notatie}}  
\newcommand{\enotat}{\end{notatie}\end{inspr}}

\newcommand{\voorb}{\begin{inspr}\begin{voorbeeld}}  
\newcommand{\evoorb}{\end{voorbeeld}\end{inspr}}

\newcommand{\stel}{\begin{inspr}\begin{stelling}}
\newcommand{\estel}{\end{stelling}\end{inspr}}

\newcommand{\prop}{\begin{inspr}\begin{propositie}}
\newcommand{\eprop}{\end{propositie}\end{inspr}}

\newcommand{\lem}{\begin{inspr}\begin{lemma}}
\newcommand{\elem}{\end{lemma}\end{inspr}}

\newcommand{\opm}{\begin{inspr}\begin{opmerking}}
\newcommand{\eopm}{\end{opmerking}\end{inspr}}

\newcommand{\voorw}{\begin{inspr}\begin{voorwaarde}}
\newcommand{\evoorw}{\end{voorwaarde}\end{inspr}}

\newcommand{\probl}{\begin{inspr}\begin{probleem}}
\newcommand{\eprobl}{\end{probleem}\end{inspr}}

\newcommand{\gev}{\begin{inspr}\begin{gevolg}}
\newcommand{\egev}{\end{gevolg}\end{inspr}}

\newcommand{\nul}{\begin{inspr}\begin{niets}}
\newcommand{\enul}{\end{niets}\end{inspr}}

\newcommand{\bew}{\vspace{-0.3cm}\begin{itemize}\item[ ] \bf Proof\rm: }
\newcommand{\ebew}{\hfill $\qed$ \end{itemize}}

\newcommand{\ssnl}{\vskip 3pt} 
\newcommand{\snl}{\vskip 7pt} 
\newcommand{\nl}{\vskip 12pt}

\newcommand{\ot}{\otimes}

\newcommand{\tl}{\triangleleft}
\newcommand{\tr}{\triangleright}

\newcommand{\tussenen}{\qquad\quad\text{and}\qquad\quad}
\newcommand{\tussen}{\qquad\quad\qquad\quad}

\newcommand{\rood}{\color{red}}
\newcommand{\blauw}{\color{blue}}

\numberwithin{thm}{section}   
\numberwithin{equation}{section} 

\newcommand{\keepcomment}[1]{}
\newcommand{\oldcomment}[1]{}

\begin{document}

\setcounter{section}{-1}  

%
%

\centerline{\bf \Large Reflections on coproducts for non-unital algebras} 
\vspace{13 pt}
\centerline{\it A.\,Van Daele \rm ($^*$)}
\vspace{30 pt}
{\bf Abstract} 
\nl 
A coproduct on a vector space $A$ is defined as a linear map $\Delta:A\to A\ot A$ satisfying coassociativity $(\Delta\ot\iota)\Delta=(\iota\ot\Delta)\Delta$. We use $\iota$ for the identity map. 
\snl
If $G$ is a {\it finite} group and if $A$ is the space of all complex functions on $G$, a coproduct on $A$ is defined by $\Delta(f)(p,q)=f(pq)$ where $p,q\in G$. We identify $A\ot A$ with complex functions on the Cartesian product $G\times G$. Coassociativity follows from the associativity of the product in $G$. 
\ssnl
Unfortunately, sometimes this  notion of a coproduct  is not the appropriate one. Just think of the above example for an {\it infinite} group. We explain this in the paper.
\ssnl
In this note, we consider the case of an algebra $A$, not necessarily unital but with a non-degenerate product. Now a coproduct is a linear map from $A$ to $M(A\ot A)$, the multiplier algebra of $A\ot A$. Unfortunately, it is no longer possible to express coassociativity in its usual form as the maps $\Delta\ot\iota$ and $\iota\ot \Delta$, defined on $A\ot A$, may no longer be defined on the multiplier algebra $M(A\ot A)$ (which in general is  bigger than the algebra $A\ot A$).
\ssnl
We will see how this problem can be overcome  in different ways. Solutions can be given so as to get various useful notions for a coproduct on a non-unital algebra. 
\ssnl
Similar problems occurs when we want to define a useful notion of a coaction in the case of non-unital algebras. We
discuss this in another paper \cite{VD-co}.
\ssnl
Not all the results we present in this paper are new. We  provide a number of references to the original papers where some of this material is treated. However, in the original papers, results are not always found in an organized form and we hope to improve that here.  Further a few solutions to some open questions are included as well as some more peculiar examples. Finally, we discuss some open problems and possible further research.
\nl 
Date: {\it 7 February  2024}
\vskip 3 cm
\hrule
\vskip 7 pt
\begin{itemize}
\item[($^*$)] Department of Mathematics, University of Leuven, Celestijnenlaan 200B,\newline
B-3001 Heverlee (Belgium). E-mail: \texttt{alfons.vandaele@kuleuven.be}
\end{itemize}
\newpage

%
%

\section{\hspace{-17pt}. Introduction} \label{s:intro}  

Recall the definition of a coalgebra. See e.g.\ Definition 2.1.3 in \cite{R}.

\defin\label{defin:coalg}
A coalgebra is a triple $(A,\Delta,\varepsilon)$ of a vector space $A$ with a coproduct $\Delta$ and a counit $\varepsilon$. The coproduct is a coassociative linear map $\Delta:A\to A\ot A$ and the counit is a linear map $\varepsilon: A\to \mathbb C$ satisfying 
\begin{equation*}
(\varepsilon\ot\iota)\Delta(a)=a
\tussenen
(\iota\ot\varepsilon)\Delta(a)=a
\end{equation*}
for all $a\in A$. 
\edefin

We use $\iota$ for the identity map. Coassociativity means that
\begin{equation}
(\Delta\ot\iota)\Delta=(\iota\ot\Delta)\Delta.\label{eqn:coass}
\end{equation}
If $A$ is not just a vector space, but an algebra, it is often assumed that $\Delta$ and $\varepsilon$ are  (unital) homomorphisms. 
However, although it is not needed for the definition, in this case, the notion is only useful for algebras with an identity. 
For a coproduct on an algebra $A$ that is not assumed to have an identity, to require that it takes values in $A\ot A$, turns out to be too restrictive. We illustrate this in Section \ref{s:copr}, see Example \ref{voorb:KG}. 
\ssnl
Instead we consider in this case linear maps from $A$ to the multiplier algebra $M(A\ot A)$. The multiplier algebra is defined when the product on $A$ is non-degenerate. This is a natural assumption on the algebra (and automatic if it has an identity).
\ssnl
But now there is a problem with coassociativity as formulated above because the maps $\Delta\ot\iota$ and $\iota\ot\Delta$ are defined on $A\ot A$, but not necessarily on the range of the coproduct in $M(A\ot A)$. Therefore  formula (\ref{eqn:coass}) as such does not make sense. 
\ssnl
The more common solution to this problem is to require, first  that the canonical maps $T_1$ and $T_2$, defined from $A\ot A$ to $M(A\ot A)$ by
\begin{equation*}
T_1(a\ot b)=\Delta(a)(1\ot b)
\tussenen
T_2(c\ot a)=(c\ot 1)\Delta(a),
\end{equation*}
have range in $A\ot A$ and then assuming that 
\begin{equation}
(c\ot 1\ot 1)(\Delta\ot\iota)(\Delta(a)(1\ot b))
=(\iota\ot\Delta)((c\ot 1)\Delta(a))(1\ot 1\ot b) \label{eqn:coass2a}
\end{equation}
for all $a,b,c\in A$. Here $1$ denotes the identity in the multiplier algebra $M(A)$. The equation (\ref{eqn:coass2a}) is the explicit form of the commutation rule $(T_2\ot\iota)(\iota\ot T_1)=(\iota\ot T_1)(T_2\ot\iota)$. It makes sense because it is assumed that these canonical maps have range in $A\ot A$.
\ssnl
Obviously, if the algebra $A$ has an identity, the condition on the range of the canonical maps is void and coassociativity as in equation (\ref{eqn:coass2a}) is equivalent with coassociativity in its usual form as in (\ref{eqn:coass}).
\ssnl
There are other possible forms of coassociativity as we explain further in this note. We will discuss the relation between these different forms.
\ssnl
In any case, it does not seem to be possible to give a suitable notion of a coproduct on a vector space, without more structure, so that it also includes the case of a coproduct on a non-unital algebra $A$ with values in $M(A\ot A)$ as above.
\ssnl
There is a similar problem when it comes to defining a coaction of a coalgebra on a vector space. This case will  be treated in a separate paper \cite{VD-co}.
\nl
In the first place, the aim of this note is to collect material that has been around for some time in papers on multiplier Hopf algebras, weak multiplier Hopf algebras and quantum hypergroups before. The focus not only lies on the notion of coassociativity for coproducts but also on the possibility to make a suitable subspace of the dual space into an algebra with the product dual to the original coproduct. All this is treated here more generally and in greater detail (and perhaps also in a more systematic way) than in these original papers. The idea is to provide  an easy and fairly complete treatment of different aspects of coassociativity as they are encountered in the theory of (weak) multiplier Hopf algebras and bialgebras. We illustrate all this with some examples.
\ssnl
From questions I got recently from young researchers who wanted to learn the subject, I have the feeling that there is a need for a note like this.
This paper is, up to a certain extend,  expository but it contains also some new results. Moreover, some new and rather special examples are include. 
\nl
\bf Content of the paper \rm
\nl
The main section is {\it Section} \ref{s:copr} where we consider various possible notions of coassociativity for coproducts, related concepts and the connection between these.
\ssnl
 In the finite-dimensional case, the concept of a coproduct is essentially the same as that of a product on the dual. However, for the coassociative maps we consider here, it is not always true that the dual space carries a product, obtained from the coproduct. One needs to consider appropriate subspaces of the dual to obtain this, on top of a few extra conditions on the original pair $(A,\Delta)$.
This is treated in {\it Section} \ref{s:dual}. 
\ssnl
 In {\it Section} \ref{s:examp} we include references to the original papers  but these cases are not treated in detail. Instead, we discuss some special cases and we give a few examples of non-regular coproducts with some indications for constructing more of such (and even more peculiar) examples. Finally,  in {\it Section} \ref{s:concl}, we reflect a little more on aspects that are still not completely understood and on the problem of finding more non-trivial examples to illustrate these aspects.
\keepcomment{
\ssnl 
\rood Check the statements about Section \ref{s:examp}.}{} 
\nl
\bf Notations and conventions, basic references \rm
\nl
We only work with (associative) algebras over the field $\mathbb C$ of complex numbers. However, it should be possible to consider other more general fields as well. The algebras need not be unital. But the product is {\it always assumed to be non-degenerate} (as a bilinear form). This is automatic if the algebra has a unit. It is also automatic if the algebra has local units. 
\ssnl
We will denote by $A'$ the space of all linear functionals on $A$.
\ssnl
Sometimes our algebras will be idempotent, i.e.\ any element is a sum of products of elements in the algebra. The condition is written as $A=A^2$. Again this is automatic if the algebra is unital, or more general, when it has local units. This condition will not be imposed on the algebra, but it is often a consequence of the other conditions that are considered.
\ssnl
We use $M(A)$ for the multiplier algebra of $A$. The multiplier algebra of a non-degenerate algebra, as we use it here, is considered in \cite{VD-mha} but it should be mentioned that it has been studied earlier (see e.g.\ \cite{Da}). We briefly recall the notion in the beginning of Section \ref{s:copr}. The identity in $M(A)$ is always denoted by $1$ while we use $\iota$ for the identity map. As a matter of fact, we will use $\iota$ for the identity map on any of the vector spaces we encounter.
\ssnl
For the notion of a coproduct as it appears in the theory of coalgebras, we refer to \cite{A}, \cite{S} and \cite{R}. For the notion of a coproduct in the theory of multiplier Hopf algebra, we refer to \cite{VD-mha} and in the setting for weak multiplier Hopf algebras, to \cite{VD-W0} and \cite{VD-W1}. 
\ssnl
The opposite algebra $A^{\text{op}}$ has the same underlying vector space as the original algebra $A$, but the product is reversed. Similarly,  the co-opposite coproduct $\Delta{^\text{cop}}$ on an algebra is obtained by flipping the original coproduct $\Delta$.
\nl
\bf Acknowledgments \rm
\nl
I would like to thank Danielle Santos Azevedo (Brasil) for drawing my attention to some aspects of coassociativity for coproducts on non-degenerate algebras that in fact I overlooked in previous work.
\ssnl
I also enjoyed my last stay in Nanjing (China) where I had the opportunity to talk about this topic and discuss this material with the master students and PhD students of my coauthor Shuanhong Wang.
\ssnl
Finally, I  had the opportonity to talk about this material in the analysis seminar at NTNU (Trondheim) during my last visit there.

%
%

 \section{\hspace{-17pt}. Coproducts and coassociativity} \label{s:copr}  
 
Let $A$ be an (associative) algebra over the field $\mathbb C$ of complex numbers. We do not assume that it is unital, but we require the product to be non-degenerate. This means that multiplication, seen as a bilinear map, is non-degenerate.
\nl
\bf The multiplier algebra of a non-degenerate algebra\rm
\nl
Recall the following definition (see the Appendix in \cite{VD-mha}).

\defin
The multiplier algebra $M(A)$ is the set of pairs $(\lambda,\rho)$ of maps from $A$ to itself satisfying $b\lambda(a)=\rho(b)a$ for all $a,b \in A$. 
\edefin

By the non-degeneracy of the product, it follows that these maps are linear and satisfy
\begin{equation}
\lambda(ac)=\lambda(a)c 
\tussenen
\rho(cb)=c\rho(b) \label{eqn:mult}
\end{equation} for all $a,b,c\in A$. 
If  $(\lambda,\rho)$ is a multiplier, then $\lambda$ is determined by $\rho$ and vice versa. 
\ssnl
There is an obvious embedding of $A$ in $M(A)$. Indeed, to an element $c$ in $A$ are associated the linear maps $\lambda:a\mapsto ca$ and $\rho:b\mapsto bc$. The pair is a multiplier because $b(ca)=(bc)a$. This is nothing else but associativity of the product on $A$. We get an {\it embedding} of $A$ in $M(A)$ precisely because the product on $A$ is assumed to be non-degenerate.

\notat
 If $x$ is a multiplier $(\lambda,\rho)$ as above, we write $xa$ for $\lambda(a)$ and $ax$ for $\rho(a)$. The defining relation then reads as coassociativity $b(xa)=(bx)a$. And the notation is consistent with the embedding of $A$ in $M(A)$.
 \enotat

It is easy to see that composition of maps yields a product on $M(A)$ and clearly the identity maps will give an identity in $M(A)$.  By definition $A$ sits in $M(A)$ as a two-sided ideal.  It is a {\it dense} ideal in the sense that for $x\in M(A)$, we have $x=0$ if $ax=0$ for all $a$ or if $xb=0$ for all $b$. In fact, $M(A)$ can be characterized as the largest unital algebra containing $A$ as a dense two-sided ideal in this way. Further it is clear that $M(A)=A$ if and only if $A$ has an identity. See \cite{VD-mha} for details.
\oldcomment{Include references - In fact, see the remark below.}{}
\ssnl
It is also possible to define the algebra $L(A)$ of left multipliers and the algebra $R(A)$ of right multipliers. Again see the Appendix in \cite{VD-mha}. Within the spirit of the notations above, a left multiplier $x$  is a linear map $a\mapsto xa$ with the property that $x(ac)=(xa)c$ for all $a,c\in A$ while a right multiplier $x$ is a linear map $b\mapsto bx$ with the property that $(cb)x=c(bx)$ for all $b,c\in A$. 
\ssnl
Having these notations in mind, it makes sense to say that a left multiplier is a multiplier if it is also a right multiplier. And similarly, a right multiplier is a multiplier if it is also a left multiplier. In other words, we think of $M(A)$ as the intersection $L(A)\cap R(A)$.
\ssnl
In \cite{VD-infmalg} we include an example of a non-degenerate, finite-dimensional algebra $A$ with the property that $L(A)$ and $R(A)$ are different from each other in the sense that there is a left multiplier that is not a multiplier and the same for a right multiplier. 
\oldcomment{
\ssnl
\rood To be verified!}{}
\nl 
Since we will be studying a coproduct on $A$, we need to consider the tensor product $A\ot A$ of $A$ with itself. The algebra $A\ot A$ is again a non-degenerate algebra and we have obvious inclusions 
\begin{equation}
A\ot A\subseteq M(A)\ot M(A)\subseteq M(A\ot A). \label{eqn:imb}
\end{equation}
The inclusions are algebra embeddings and in general, except if $A$ has an identity, the two inclusions are strict. See e.g.\ Example \ref{voorb:KG} below.

\opm
The notion of a multiplier in the case of a non-unital algebra 
with a non-degenerate product as above 
was introduced in the first paper on multiplier Hopf algebras (\cite{VD-mha}). The treatment was inspired by the concept as used in the theory of operator algebras (see e.g. Section 3.12 in \cite{P}). I have a background in that field. However, as it was pointed out to me later, the notion in fact had been introduced already before in 1969  by J.\ Dauns in \cite{Da}.
\eopm

We also recall the notion in detail in \cite{VD-lumalg} where some examples are found. In another paper,  \cite{VD-infmalg}, where we treat infinite matrix algebras, also the multiplier algebra is obtained.
\ssnl
\oldcomment{
\ssnl
\rood We have to check the above references to these papers.}{}
\keepcomment{
\ssnl
For the moment, we leave this introduction on multiplier Hopf algebras here. One can also find it in  \cite{VD-lumalg}. There you find more details. Wait until we know more about these two papers by me and Joost. \rood In het oog houden!}{}
\nl
\bf The notion of a coproduct on a non-degenerate algebra \rm 
\nl
By a coproduct $\Delta$ on a vector space $A$  we usually mean a linear map from $A$ to $A\ot A$ satisfying coassociativity $(\Delta\ot\iota)\Delta=(\iota\ot\Delta)\Delta$. As mentioned before, $\iota$ is the identity map from $A$ to itself. The maps $\Delta\ot\iota$ and $ \iota\ot\Delta$ are the obvious maps from $A\ot A$ to $A\ot A\ot A$. We consider composition of maps. 
\oldcomment{Wat opletten met herhalingen!}{}
\ssnl
In principle, we could use this definition for a coproduct on a non-degenerate algebra $A$ without a unit. However, this  gives a far too restrictive concept as the following motivating example indicates.

\voorb \label{voorb:KG}
Let $G$ be a group and assume that it is {\it not finite}. Let $A$ be the space $K(G)$  of complex functions with finite support in $G$. It is a non-degenerate algebra (for the pointwise product) and it will not have a unit. The multiplier algebra $M(A)$ is canonically identified with the algebra of all complex functions (again with pointwise product). The tensor product algebra $A\ot A$ is identified with the algebra $K(G\times G)$ of complex functions with finite support on the Cartesian product $G\times G$ and its multiplier algebra $M(A\ot A)$ is then the algebra of all complex functions on $G\times G$. 
\ssnl
Remark in passing that for this example, the embeddings in (\ref{eqn:imb}) are indeed strict.
\ssnl
The multiplication in $G$ induces a homomorphism from $A$ to $M(A\ot A)$ by the formula
\begin{equation*}
\Delta(f)(p,q)=f(pq)
\end{equation*}
where $f\in K(G)$ and $p,q\in G$. The function $\Delta(f)$ will never be in $A\ot A$, except when $f$ is $0$. Indeed, suppose that $r\in G$ and that $ f(r)\neq 0$. Then 
\begin{equation*}
\Delta(f)(p,p^{-1}r)=f(pp^{-1}r)=f(r)
\end{equation*} 
and this is non-zero for all $p$. As $G$ is supposed to be infinite, the function $\Delta(f)$ will not have finite support. Hence $\Delta(A)$ is not contained in $A\ot A$.
\ssnl
It would not help to take all complex functions for $A$ as in general, it will not be possible to write $\Delta(f)=\sum f_i\ot g_i$ with a finite sum.
\evoorb

If $G$ is not a group, but just a set with an associative multiplication, the situation might be different. Take e.g.\ for $G$ the natural numbers $\mathbb N$ with addition. If we take for $A$ again the algebra of functions with finite support, we will have $\Delta(A)\subseteq A\ot A$ (with $\Delta$ as defined above). Still, this will not be true when we take for $A$ the algebra of all complex functions.
\ssnl
If $G$ is a finite group, the above problem obviously does not occur. But if we stick to the usual definition of a coproduct, we can not consider the case where $G$ is infinite. this would be very restrictive from various viewpoints.

\opm
There are cases where the problem can be avoided. Take e.g.\ a matrix group $G$ and  for $A$ the subalgebra of the algebra $C(G)$ of all complex functions, spanned by the matrix elements (seen as functions on the group). With the coproduct defined as before, we do have $\Delta(A)\subseteq A\ot A$. This choice has other disadvantages. If $G$ is not a {\it compact} matrix group, the matrix elements will not give rise to integrable functions for the Haar measure and consequently, the algebra $A$ will not have integrals. This is one of the reasons why it is not satisfactory, also for such a group, to restrict to a coproduct with values in the tensor product.
\eopm

This leads to the following problem.

\probl Let $A$ be a non-degenerate algebra and $\Delta$ a linear map from $A$ to $M(A\ot A)$. In general, the linear maps $\Delta\ot\iota$ and $\iota\ot \Delta$ can be defined on $A\ot A$, but there are no natural and canonical extensions of these maps to $M(A\ot A)$. Therefore coassociativity, in its usual form $(\Delta\ot\iota)\Delta=(\iota\ot\Delta)\Delta$, has no meaning.
\ssnl
We need some extra regularity conditions in order to be able to formulate coassociativity. 
\eprobl

There are several possibilities. Some of them are of a very different nature, within different situations.
\nl
\bf Regularity of the canonical maps \rm
\nl
A first set of possible conditions that allows the notion of coassociativity uses the following definition for a {\it linear map} $\Delta:A\to M(A\ot A)$.

\defin\label{defin:1.2}
Let $A$ be a non-degenerate algebra and $\Delta$ a linear map from $A$ to\\
$M(A\ot A)$. We associate four linear maps $T_1$, $T_2$, $T_3$ and $T_4$ from $A\ot A$ to $M(A\ot A)$. We define
\begin{align*}
T_1(a\ot b)&=\Delta(a)(1\ot b)
\tussenen 
T_2(c\ot a)=(c\ot 1)\Delta(a) \\
T_3(a\ot b)&=(1\ot b)\Delta(a)
\tussenen
T_4(c\ot a)=\Delta(a)(c\ot 1)
\end{align*}
for all $a,b,c\in A$. We call these maps the {\it canonical maps} (associated with $\Delta$). Such a canonical map is called {\it regular} if it has range in $A\ot A$.
\edefin
Here $1$ is the identity in the multiplier algebra $M(A)$. We use that elements of the form $1\ot b$ and $c\ot 1$ belong to $M(A\ot A)$ via the embeddings of $A$ in $M(A)$ and of $M(A)\ot M(A)$ in $M(A\ot A)$.
\ssnl
Regularity of the canonical maps is natural as we can see in the case of the group example in Example \ref{voorb:KG}. 

\voorb
Indeed, for the map $T_1$ we find $T_1(f)(p,q)=f(pq,q)$ when $f$ is a function with finite support in $G\times G$. For the result to be non zero, $q$ has to lie in a finite set and as this is also true for $pq$, by the group property, we have that moreover $p$ is forced to lie in a finite set. Hence $T_1(f)$ is again a function with finite support. So $T_1$ is regular. Similarly for the three other canonical maps.
\evoorb

For this example, the four canonical maps are regular. In Section \ref{s:examp}, we give some examples where only some of the canonical maps are regular, see e.g.\ Example \ref{voorb:4.1}, \ref{voorb:4.3}, \ref{voorb:4.4} and others in Section \ref{s:examp}. 
\ssnl
This leads  to the first possible definitions of coassociativity.

\defin\label{defin:coass} 
i) If the maps $T_1$ and $T_2$ are regular, we say that $\Delta$ is coassociative if 
\begin{equation}
(c\ot 1\ot 1)(\Delta\ot\iota)(\Delta(a)(1\ot b))
=(\iota\ot\Delta)((c\ot 1)\Delta(a))(1\ot 1\ot b)\label{eqn:1.3}
\end{equation}
for all $a,b,c\in A$. 
\ssnl
ii) If the maps $T_3$ and $T_4$ are regular, we say that $\Delta$ is coassociative if
\begin{equation}
(\Delta\ot\iota)((1\ot b)\Delta(a))(c\ot 1\ot 1)
=(1\ot 1\ot b)(\iota\ot\Delta)(\Delta(a)(c\ot 1))\label{eqn:1.3a}
\end{equation}
for all $a,b,c\in A$. 
\edefin

Remark that regularity of the maps $T_1$ and $T_2$ is needed  for Equation (\ref{eqn:1.3}) to make sense. Similarly we need the regularity of $T_3$ and $T_4$ for Equation \ref{eqn:1.3a}. One form is equivalent with the other when we replace $A$ by $A^\text{op}$ or $\Delta$ by $\Delta^\text{cop}$.
\ssnl
Equation (\ref{eqn:1.3}) can be reformulated as 
$(T_2\ot\iota)(\iota\ot T_1)=(\iota\ot T_1)(T_2\ot\iota)$. Similar alternative expressions are possible for the other cases. 
\ssnl
With these conditions, we get the  first definition of a coproduct on a non-degenerate algebra.

\defin\label{defin:copr}
Let $A$ be a non-degenerate algebra and $\Delta$ a linear map from $A$ to $M(A\ot A)$. We call it a {\it coproduct} (or {\it comultiplication}) if $T_1$ and $T_2$ are regular and if $(T_2\ot\iota)(\iota\ot T_1)=(\iota\ot T_1)(T_2\ot\iota)$ holds. 
\edefin

We encounter this notion as in item i) of Definition \ref{defin:copr} for the first time in \cite{VD-mha} where multiplier Hopf algebras are introduced. Later we use it for algebraic quantum hypergroups (see e.g.\ \cite{De-VD}), weak multiplier Hopf algebras (see e.g.\ \cite{VD-W1}) and in many other situations. The condition as in item ii) of Definition \ref{defin:copr} is illustrated in some of our examples in Section \ref{s:examp}, see Proposition  \ref{prop:3.9d}.
\ssnl
There are still other possibilities. Consider the following two cases. As before, $A$ is a non-degenerate algebra and $\Delta:A\to M(A\ot A)$  a linear map.

\defin\label{defin:1.11a} 
i) Assume that $T_1$ and $T_4$ are regular. Then we call $\Delta$ coassociative if
\begin{equation*}
((\Delta\ot\iota)(\Delta(a)(1\ot b)))(c\ot 1\ot 1)
=((\iota\ot\Delta)(\Delta(a)(c\ot 1)))(1\ot 1\ot b)
\end{equation*}
for all $a,b,c\in A$.
\ssnl
ii) \oldcomment{This is wrong!}{}
Assume that $T_2$ and $T_3$ are regular. Then we call $\Delta$ coassociative if
\begin{equation*}
(c\ot 1\ot 1)((\Delta\ot\iota)((1\ot b)\Delta(a)))
=(1\ot 1\ot b)((\iota\ot\Delta)((c\ot 1)\Delta(a)))
\end{equation*}
for all $a,b,c\in A$.
\edefin

These two forms of coassociativity are used in the theory of multiplier Hopf algebroids, see \cite{T-VD}.
\ssnl
It is also possible to use the explicit forms of the canonical maps, but the formulas are not so nice. Moreover, whereas the expressions for $T_1$ and $T_2$ are more or less standard, the choices for the two other ones, $T_3$ and $T_4$ are not.
\ssnl
There are still other combinations of two canonical maps possible but we will not consider them further.
\ssnl
In the following case, all these forms are equivalent with each other.

\defin Let $A$ be a non-degenerate algebra and $\Delta$ a linear map from $A$ to $M(A\ot A)$. We call $\Delta$ {\it regular} if all four canonical maps are regular. 
\edefin

Indeed, using the non-degeneracy of the product, when all canonical maps are regular, one can easily show that $\Delta$ is coassociative in any of the two forms in Definition \ref{defin:coass} if and only if it is coassociative in the sense of any of the two possibilities in  Definition \ref{defin:1.11a}.
\nl
\bf Coassociativity with only one regular canonical map \rm
\nl
 Under certain other conditions, it is sufficient that only one canonical map is regular. We treat this in the following definition for the case where $T_1$ is regular. Observe that we now require that $\Delta:A\to M(A\ot A)$ is a \emph{homomorphism} and not just a linear map.

\defin\label{defin:1.6}
Let $A$ be a non-degenerate algebra and $\Delta:A\to M(A\ot A)$ a homomorphism. Assume that the map $T_1$ is regular. We call $\Delta$ coassociative if
\begin{equation}
((\iota\ot\Delta)(\Delta(a)(1\ot p)))(1\ot 1\ot q)
=\sum_i(\Delta\ot\iota)(\Delta(a)(1\ot c_i))(1\ot b_i\ot 1)\label{eqn:1.7a}
\end{equation}
where $\Delta(p)(1\ot q)=\sum_i b_i\ot c_i$. 
\edefin

We can show the following equivalence.

\prop
Assume that $\Delta$ is a homomorphism and that $T_1$ and $T_2$ are both regular. If $\Delta$ is coassociative in the sense of  Definition \ref{defin:coass} then it is also coassociative in the sense of Definition \ref{defin:1.6}. If moreover $T_1$ is surjective, also the converse is true.
\eprop

\bew
i) Assume that $\Delta$ is coassociative in the sense of Definition \ref{defin:coass}. Let $a,b,p,q\in A$. Then we have
\begin{align*}
(c\ot 1\ot 1)& ((\iota\ot\Delta)(\Delta(a)(1\ot p)))(1\ot 1\ot q)\\
&= (\iota\ot\Delta)((c\ot 1)\Delta(a)(1\ot p))(1\ot 1\ot q)\\
&=(\iota\ot\Delta)((c\ot 1)\Delta(a))(1\ot (\Delta(p)(1\ot q))).
\end{align*}
We have used that $\Delta$ is a homomorphism. Now write $\Delta(p)(1\ot q)=\sum_i b_i\ot c_i$. Then we get
\begin{align*}
(c\ot 1\ot 1)& ((\iota\ot\Delta)(\Delta(a)(1\ot p)))(1\ot 1\ot q)\\
&=\sum_i (\iota\ot\Delta)((c\ot 1)\Delta(a))(1\ot b_i \ot c_i)\\
&=\sum_i (\iota\ot\Delta)((c\ot 1)\Delta(a))(1\ot 1 \ot c_i)(1\ot b_i\ot 1)\\
&= \sum_i(c\ot 1\ot 1)((\Delta\ot \iota)(\Delta(a)(1 \ot c_i)))(1\ot b_i\ot 1).
\end{align*}
Now we have used coassociativity as in Equation \ref{eqn:1.3} of Definition \ref{defin:coass}. We can cancel $(c\ot 1\ot 1)$ and we get 
\begin{equation*}
((\iota\ot\Delta)(\Delta(a)(1\ot p)))(1\ot 1\ot q)
= \sum_i(\Delta\ot \iota)(\Delta(a)(1 \ot c_i))(1\ot b_i\ot 1).
\end{equation*}
We see that $\Delta$ is coassociative in the sense of Definition \ref{defin:1.6}.
\ssnl
ii) To prove the converse let $a,b,c, ,p,q\in A$ and $\Delta(p)(1\ot q)=\sum b_i\ot c_i$. Assuming coassociativity as in Definition \ref{defin:1.6}, from the calculation in the previous item we now find that 
\begin{align*}
\sum_i (\iota\ot\Delta)&((c\ot 1)\Delta(a))(1\ot 1 \ot c_i)(1\ot b_i\ot 1)\\
&= \sum_i(c\ot 1\ot 1)((\Delta\ot \iota)(\Delta(a)(1 \ot c_i)))(1\ot b_i\ot 1).
\end{align*}
By the assumption that $T_1$ is surjective, we can replace $b_i, c_i$ by any pair $b,d$ and we obtain
\begin{align*}
(\iota\ot\Delta)&((c\ot 1)\Delta(a))(1\ot 1 \ot d)(1\ot b\ot 1)\\
&= (c\ot 1\ot 1)(\Delta\ot \iota)(\Delta(a)(1 \ot d))(1\ot b\ot 1).
\end{align*}
Now we can cancel $b$ and get
\begin{equation*}
(\iota\ot\Delta)((c\ot 1)\Delta(a))(1\ot 1 \ot d)
= (c\ot 1\ot 1)(\Delta\ot \iota)(\Delta(a)(1 \ot d)).
\end{equation*}
This completes the proof.
\ebew

In the next section, we will formulate another notion of coassociativity that, in a way, incorporates all these notions involving the canonical maps. We will do this in connection with the construction of a dual algebra $B$ where the product is obtained from the coproduct on $A$, see Definition \ref{defin:2.8}. \oldcomment{\rood Include detailed references}

\nl
\bf Full coproducts and counits \rm
\nl
Up to now, we did not formulate any condition on the coproduct that prevents it to be completely trivial. A first possibility is the existence of a counit.

\defin\label{defin:1.15c} 
Suppose that the maps $T_1$ and $T_2$ are regular. Then one defines a counit as a linear map $\varepsilon:A\mapsto \mathbb C$ satisfying
\begin{equation*}
(\varepsilon\ot\iota)(\Delta(a)(1\ot c))=ac
\tussenen
(\iota\ot\varepsilon)((c\ot 1)\Delta(a))=ca
\end{equation*}
for all $a,c\in A$.
\edefin

Under certain natural conditions, such a counit is unique. This is easily seen to be the case if $\varepsilon$ is a homomorphism. But it is also true if $\Delta$ is full (see Definition \ref{defin:1.16c} below). \oldcomment{Add a concrete reference}{}
\ssnl
If e.g.\ only $T_1$ is regular, only the first equality makes sense. We could speak of a left counit in that case. See examples in Section \ref{s:examp} later. \keepcomment{\rood Add reference. More comments?}{}

\defin\label{defin:1.16c}
Assume that $T_1$ and $T_2$ are regular. Then the coproduct is called full if the smallest subspaces $V$ and $W$ of $A$, satisfying
\begin{equation*}
\Delta(a)(1\ot c)\in V\ot A
\tussenen
(c\ot 1)\Delta(a)\in A\ot W
\end{equation*}
for all $a,c\in A$, are actually $A$ itself.
\edefin

It is not hard to show that the coproduct is full if and only if the linear span of elements of the form $(f\ot\iota)((c\ot 1)\Delta(a))$ with $f\in A'$ and $a,c\in A$, as well as the linear span of elements of the form $(\iota\ot f)(\Delta(a)(1\ot c))$ are all of $A$.
\ssnl
If only $T_1$ is regular, we can speak about the left leg of $\Delta$. If only $T_2$ is regular, we can define the right leg. 
\ssnl
For more information about the notion of fullness of a coproduct and the relation with the counit, see e.g.\ Section 1 of \cite{VD-W0}
\ssnl 
In the following item, we have a notion of coassociativity, not involving the canonical maps.  It is another way to guarantee that $\Delta$ is non-trivial.
\nl
\bf Coassociativity for non-degenerate coproducts \rm
\nl
We will now treat coassociativity for {\it non-degenerate coproducts}. This notion is  of a different type and does not even involve regularity of the canonical maps. On the other hand, again it only makes sense when $\Delta$ is a homomorphism.

\defin\label{defin:nd}
Let $A$ be a non-degenerate algebra and $\Delta$ a homomorphism from $A$ to $M(A\ot A)$. We call $\Delta$ {\it non-degenerate} if 
\begin{equation*}
\Delta(A)(A\ot A)=A\ot A
\tussenen 
(A\ot A)\Delta(A)=A\ot A.
\end{equation*}
\edefin

Because $\Delta$ maps $A$ into $M(A\ot A)$, the spaces $\Delta(A)(A\ot A)$ and $(A\ot A)\Delta(A)$ are subspaces of $A\ot A$. The requirement is that they are actually all of $A\ot A$. This means that for any pair $a,b\in A$ there exists a finite number of element $c_i,a_i,b_i$ in $A$ such that $a\ot b=\sum_i\Delta(c_i)(a_i\ot b_i)$ and similarly for the equality $(A\ot A)\Delta(A)=A\ot A$.
\ssnl
We will show later that the coproduct in a multiplier Hopf algebra is non-degenerate, see Proposition \ref{prop:1.27d}.

\prop\label{prop:ext}
If $\Delta:A\to M(A\ot A)$ is a non-degenerate homomorphism, then there is a unique homomorphism $\Delta_1:M(A)\to M(A\ot A)$ extending $\Delta$ on $A$. The extension $\Delta_1$ is unital.
\eprop

\bew
Let $m\in M(A)$. 
\vskip 2pt
i) We claim that we can define a linear map from $\lambda$ from $A\ot A$ to itself  by 
\begin{equation*}
\lambda(\Delta(c)(a\ot b))=\Delta(mc)(a\ot b)
\end{equation*}
 for all $a,b,c\in A$. To prove this, assume that we have a finite number of elements $a_i,b_i,c_i$ in $A$ such that
$\sum_i\Delta(c_i)(a_i\ot b_i)=0$. Take any $p,q,r\in A$. Then 
\begin{align*}
(p\ot q)\Delta(r)\sum_i\Delta(mc_i)(a_i\ot b_i)
&=\sum_i(p\ot q)\Delta(rmc_i)(a_i\ot b_i)\\
&=(p\ot q)\Delta(rm)\sum_i\Delta(c_i)(a_i\ot b_i)=0.
\end{align*}
Because it is assumed that $(A\ot A)\Delta(A)=A\ot A$, we get that also 
\begin{equation*}
x\sum_i\Delta(mc_i)(a_i\ot b_i)=0
\end{equation*}
 for all $x\in A\ot A$. By the non-degeneracy of the product in $A\ot A$, it follows that $\sum_i\Delta(mc_i)(a_i\ot b_i)=0$. Then the claim follows. Because $\Delta(A)(A\ot A)=A\ot A$, the map $\lambda$ is defined on all of $A\ot A$.
\vskip 2pt
ii) Similarly, there is a map $\rho:A\ot A\to A\ot A$ defined by
\begin{equation*}
\rho((a\ot b)\Delta(c))=(a\ot b)\Delta(cm)
\end{equation*}
 for all $a,b,c\in A$.
 \vskip 2pt
 iii) We have 
 \begin{align*}
(p\ot q)\Delta(r)\lambda(\Delta(c)(a\ot b))
&=(p\ot q)\Delta(r)\Delta(mc)(a\ot b)\\
&=(p\ot q)\Delta(rmc)(a\ot b)\\
&=(p\ot q)\Delta(rm)\Delta(c)(a\ot b)\\
&=\rho((p\ot q)\Delta(r))(\Delta(c)(a\ot b))
\end{align*}
for all $a,b,c$ and $p,q,r$ in $A$. This proves that the pair $(\lambda,\rho)$ is a multiplier in $M(A\ot A)$. We denote it by $\Delta_1(m)$.
\vskip 2pt
It is now more or less obvious to show that $\Delta_1$ is a homomorphism from $M(A)$ to $M(A\ot A)$, that it extends $\Delta$ on $A$ and that it is the unique homomorphism with this property. It is also clear that $\Delta_1(1)=1\ot 1$ in $M(A\ot A)$.
\ebew

Remark that here we do need that $\Delta$ is a homomorphism and not just a linear map. This was sufficient to formulate coassociativity using regularity of the canonical maps but it will not be sufficient here.

\opm\label{opm:ext}
The above result is easily seen to be a special case of a more general result that says that a non-degenerate homomorphism $\gamma:A\to M(B)$,  where $A$ and $B$ are non-degenerate algebras, has a unique extension to a unital homomorphism from $M(A)$ to $M(B)$. See e.g.\ Proposition A.5 in \cite{VD-mha}.
\eopm

It is a common practice to denote this extension $\Delta_1$ again by the symbol $\Delta$ used for the original homomorphism. We will do this in what follows.

\opm\label{opm:1.20}
Unfortunately, there is a slight problem with this notion and the corresponding result. Consider e.g.\ the embedding  from $A$ in $M(A)$. This is a homomorphism and it is non-degenerate in the sense of Definition \ref{defin:nd} if and only if $AA=A$, that is when $A$ is idempotent. But of course, also if $A$ is not idempotent, the embedding extends to the identity map from $M(A)$ to itself. Fortunately, in all the relevant cases, we are working with idempotent algebras. See a result proven further in this item, Proposition \ref{prop:idempotent}. See also a remark following Proposition \ref{prop:3.3d} in Section \ref{s:examp}.
\eopm

Some more reflections on this problem can be found in \cite{VD-lumalg}. 
\oldcomment{\rood Provide a more specific reference.}{} 
\ssnl
A similar problem occurs when we want to treat coassociativity in this context as we see in what follows.
\ssnl
Let $\Delta$ be a non-degenerate homomorphism from $A$ to $M(A\ot A)$. Consider the homomorphisms $\Delta\ot \iota$ and $\iota\ot\Delta$ from $A\ot A\to M(A\ot A\ot A)$. If the algebra $A$ is idempotent, these homomorphisms are still non-degenerate. The general procedure formulated above, can be applied and we get unique extensions, still denoted by the same symbols, to unital homomorphisms from $M(A\ot A)$ to $M(A\ot A\ot A)$.
\ssnl
However,  it is not necessary to know that $A$ is idempotent to obtain these extensions. We show this in the next proposition.

\prop \label{prop:ext2}
Assume that $\Delta:A\to M(A\ot A)$ is a non-degenerate homomorphism. Then the homomorphisms $\Delta\ot \iota$ and $\iota\ot\Delta$ from $A\ot A\to M(A\ot A\ot A)$ have unique extensions to unital homomorphisms from $M(A\ot A)$ to $M(A\ot A\ot A)$. 
\eprop

\bew
The proof is very much as the proof of Proposition \ref{prop:ext}.
\vskip 2pt
Consider $m\in M(A\ot A)$. 
We define  $(\Delta\ot \iota)(m)$ as a left multiplier of $A\ot A\ot A$ by the formula
\begin{equation*}
((\Delta\ot\iota)(m))(\Delta(c)(a\ot b)\ot d)=(\Delta\ot\iota)(m(c\ot d))(a\ot b\ot 1)
\end{equation*}
and as a right multiplier by the formula
\begin{equation*}
((a\ot b)\Delta(c)\ot d)(\Delta\ot\iota)(m))=(a\ot b\ot 1)(\Delta\ot\iota)((c\ot d)m)
\end{equation*}
where $a,b,c,d\in A$. Then we continue as in the proof of Proposition \ref{prop:ext}.
\ebew
 
Now we can formulate another notion of coassociativity.

\defin\label{defin:coass2}
Let $\Delta$ be a non-degenerate homomorphism from $A$ to $M(A\ot A)$. We call $\Delta$ coassociative if $(\Delta\ot\iota)\Delta=(\iota\ot\Delta)\Delta$. 
\edefin

Here $\Delta\ot\iota$ and $\iota\ot\Delta$ are the extensions obtained in Proposition \ref{prop:ext2}. 
\ssnl
The formula is required to hold on $A$. But it can be shown that it will also hold on $M(A)$ when also the extension of $\Delta$ to $M(A)$ is used as obtained in Proposition \ref{prop:ext}. This is done in the next proposition.

\prop
Let $\Delta:A\to M(A\ot A)$ be a non-degenerate homomorphism. If $(\Delta\ot\iota)\Delta(a)=(\iota\ot\Delta)\Delta(a)$ holds for all $a\in A$ then also $(\Delta\ot\iota)\Delta(m)=(\iota\ot\Delta)\Delta(m)$ for all $m\in M(A)$.
\eprop

\bew
Let $m$ be in $M(A)$ and $a$ in $A$. Then
\begin{align*}
(\Delta\ot\iota)\Delta(m)(\Delta\ot\iota)\Delta(a)
&=(\Delta\ot\iota)(\Delta(m)\Delta(a)) \\
&=(\Delta\ot\iota)\Delta(ma).
\end{align*}
We have used that $\Delta\ot\iota$ is a homomorphism on $M(A\ot A)$ and that $\Delta$ is a homomorphism on $M(A)$. We have a similar result for $(\iota\ot\Delta)\Delta(m)$ and from coassociativity, we find that 
\begin{align*}
(\Delta\ot\iota)\Delta(m)(\Delta\ot\iota)\Delta(a)
&=(\iota\ot\Delta)\Delta(m)(\iota\ot\Delta)\Delta(a)\\
&=(\iota\ot\Delta)\Delta(m)(\Delta\ot\iota)\Delta(a).
\end{align*}
Let us now denote $(\Delta\ot\iota)\Delta(m)$ by $M$ and $(\iota\ot\Delta)\Delta(m)$ by $N$. If we multiply the previous equation by $(\Delta\ot\iota)(b\ot c)$ and use that $\Delta(A)(A\ot A)=A\ot A$, we find that
\begin{equation*}
M(\Delta\ot\iota)(p\ot q)=N(\Delta\ot\iota)(p\ot q)
\end{equation*}
for all $p,q\in A$. Next we multiply with $r\ot s\ot 1$ and use again that  $\Delta(A)(A\ot A)=A\ot A$. It follows that 
\begin{equation*}
M(u\ot v\ot q)=N(u\ot v\ot q)
\end{equation*}
for all $u,v,q\in A$. Hence $M=N$ and this completes the proof.
\ebew

We now show that the previous {\it notions of coassociativity}, like the first one as in Definition \ref{defin:coass} and the second one as in Definition \ref{defin:1.11a}, are {\it compatible} with Definition \ref{defin:coass2}. Here is a more precise formulation.

\prop\label{prop:equiv}
Let $\Delta:A\to M(A\ot A)$ be a non-degenerate homomorphism. Assume that the maps $T_1$ and $T_2$ are regular or that $T_3$ and $T_4$ are regular. Then $\Delta$ is coassociative in the sense of Definition \ref{defin:coass2} if and only if it is coassociative in the sense of Definition \ref{defin:coass}. A similar result holds if any of the other pairs of canonical maps, $T_1$ and $T_4$ or $T_2$ and $T_3$ are regular (as in Definition \ref{defin:1.11a}).
\eprop

\bew
i) Start with a non-degenerate homomorphism $\Delta:A\to M(A\ot A)$ satisfying coassociativity $(\Delta\ot\iota)\Delta=(\iota\ot\Delta)\Delta$ as in Definition \ref{defin:coass2}. Take $a,b,c\in A$. We get 
\begin{align*}
(c\ot 1\ot 1)((\Delta\ot\iota)\Delta(a))&(1\ot 1\ot b)= \\
&(c\ot 1\ot 1)((\iota\ot\Delta)\Delta(a))(1\ot 1\ot b).
\end{align*}
We claim that 
\begin{align*}
&(c\ot 1\ot 1)((\iota\ot\Delta)(m))=(\iota\ot\Delta)((c\ot 1)m)\\
&((\Delta\ot\iota)(m))(1\ot 1\ot b)=(\Delta\ot\iota)(m(1\ot b))
\end{align*}
for all $m\in M(A\ot A)$. This will imply that
\begin{equation}
(c\ot 1\ot 1)(\Delta\ot\iota)(\Delta(a)(1\ot b))
=(\iota\ot\Delta)((c\ot 1)\Delta(a))(1\ot 1\ot b)\label{eqn:coass5}
\end{equation}
and so we have coassociativity of $\Delta$ as in Definition \ref{defin:coass}.
\vskip 3pt
To prove the claim, we multiply with $1\ot\Delta(a)$ from the right in the case of the first formula and from the left in the case of the second formula and we use that the formulas hold for $m$ in $A$.
\vskip 4pt
ii) 
Conversely, suppose that for all $a,b,c$ in $A$ we have the above equality \ref{eqn:coass5}. Essentially by the same arguments as above, it will follow that 
\begin{equation*}
(c\ot 1\ot 1)((\Delta\ot\iota)\Delta(a))(1\ot 1\ot b)
=(c\ot 1\ot 1)((\iota\ot\Delta)\Delta(a))(1\ot 1\ot b).
\end{equation*}
This holds for all $b,c$ and so
\begin{equation*}
(\Delta\ot\iota)\Delta(a)
=(\iota\ot\Delta)\Delta(a).
\end{equation*} 
This proves the other direction.
\ssnl
In a similar way, we can prove the result for the other cases.
\ebew
A similar result is true if we have coassociativity with only one of regular canonical map (as in Definition \ref{defin:1.6}). 
\ssnl
Before we pass to a weak condition of non-degeneracy for a coproduct, we would like to add the following result and a remark about the case of a multiplier Hopf algebra.

\prop\label{prop:idempotent}
 For a non-degenerate homomorphism $\Delta:A\to M(A\ot A)$, if any of the canonical maps is regular, then $A$ has to be idempotent. 
\eprop
\bew
Indeed, suppose e.g.\ that $T_1$ is regular. Let $\omega$ be a linear functional on $A$ that is zero on $A^2$. Then $(\omega\ot\iota)(\Delta(a)(b\ot c))=0$ for all $a,b,c\in A$ because $\Delta(a)(1\ot c)\in A\ot A$. Now use that elements of the form $\Delta(a)(b\ot c)$ span all of $A\ot A$. It follows that $(\omega\ot\iota)(p\ot q)=0$ for all $p,q\in A$. Hence $\omega=0$. This proves that $A^2=A$. A similar argument works when any of the four canonical maps is regular.
\ebew

This means that in this situation, we can obtain the extensions $\Delta\ot\iota$ and $\iota\ot\Delta$ simply by the general procedure as mentioned in a remark before Proposition \ref{prop:ext2}. 
\ssnl
We have variations of this result.

\prop
Assume that $T_1$ and $T_2$ are regular. If any of these canonical maps has range all of $A\ot A$, then $A$ is idempotent.
\eprop

\bew
Let $\omega$ be a linear functional on $A$ that is zero on $A^2$. Then
\begin{equation*}
(\iota\ot\omega)((c\ot 1)\Delta(a)(1\ot b))=0
\end{equation*}
for all $a,b,c\in A$. We use that $(c\ot 1)\Delta(a)\in A\ot A$ so that $(c\ot 1)\Delta(a)(1\ot b)\in A\ot A^2$. However, because $\Delta(a)(1\ot b)\in A\ot A$ we can cancel $c$ and get that also
\begin{equation*}
(\iota\ot\omega)(\Delta(a)(1\ot b))=0.
\end{equation*}
If $T_1$ has range equal to $A\ot A$ and we see that $\omega$ is zero on all of $A$. Therefore, $A=A^2$. 
\ssnl
A similar argument works when the range of $T_2$ is $A\ot A$.
\ebew

\prop\label{prop:1.27d}
If the two canonical maps $T_1$ and $T_2$ are regular and have range $A\ot A$, then $\Delta$ is non-degenerate. In particular, when $(A,\Delta)$ is a multiplier Hopf algebra, the coproduct is non-degenerate.
\eprop
\bew
We know that $A$ has to be idempotent. And if the maps $T_1$ and $T_2$ have range $A\ot A$, we have 
\begin{align*}
\Delta(A)(A\ot A)&=\Delta(A)(1\ot A)(A\ot 1)=(A\ot A)(A\ot 1)=A^2\ot A=A\ot A\\
(A\ot A)\Delta(A)&=(1\ot A)(A\ot 1)\Delta(A)=(1\ot A)(A\ot A)=A\ot A^2=A\ot A.
\end{align*}
When $(A,\Delta)$ is a multiplier Hopf algebra, the canonical maps $T_1$ and $T_2$ are bijective maps from $A\ot A$ to itself. Then the previous result applies. 
\ebew

This result is not found in the literature on multiplier Hopf algebras, although it may have been implicitly used. \keepcomment{\rood More comments?}{}

\nl
\bf Weakly non-degenerate coproducts \rm
\nl
In the previous item, we considered coassociativity of a coproduct if the coproduct is assumed to be non-degenerate. This is the situation that occurs in the study of multiplier Hopf algebras as we have seen in Proposition  \ref{prop:1.27d}. In the case of weak multiplier Hopf algebras however, we do not work with non-degenerate coproducts but with a weaker condition. We recall it below.

\defin\label{defin:wnd}
Let $A$ be a non-degenerate algebra and $\Delta$ a homomorphism from $A$ to $M(A\ot A)$. We say that $\Delta$ is {\it weakly non-degenerate} if there is an idempotent element $E$ in $M(A\ot A)$ so that 
\begin{equation*}
\Delta(A)(A\ot A)=E(A\ot A)
\tussenen
(A\ot A)\Delta(A)=(A\ot A)E.
\end{equation*}
\edefin

If such an idempotent exists, it is unique. It satisfies $E\Delta(a)=\Delta(a)$ and $\Delta(a)E=\Delta(a)$ for all $a\in A$. Moreover, it is the smallest idempotent with this property. By this we mean that, if $F$ is another idempotent in $M(A\ot A)$ so that $F\Delta(a)=\Delta(a)$ and $\Delta(a)F=\Delta(a)$ for all $a\in A$, then $EF=E$ and $FE=E$. See \cite{VD-W1}.
\ssnl
If $E=1$, then $\Delta$ is non-degenerate and it has a unique extension to a homomorphism from $M(A)$ to $M(A\ot A)$ (as we proved in Proposition \ref{prop:ext}). We will now show that this result can be proven, also under this weaker condition. However, now the extension will no longer be unital. 

\prop\label{prop:ext3}
If $\Delta:A\to M(A\ot A)$ is a weakly non-degenerate homomorphism, then there is a unique homomorphism $\Delta_1:M(A)\to M(A\ot A)$ extending $\Delta$ on $A$, satisfying $\Delta_1(1)=E$. 
\eprop

\bew
The proof is not very different from the one in the more restrictive situation.
\ssnl
Take $m\in M(A)$. Define linear maps $\lambda$ and $\rho$ from $E(A\ot A)$, respectively $(A\ot A)E$, to $A\ot A$ by 
\begin{align}
\lambda(\Delta(c)(a\ot b))&=\Delta(mc)(a\ot b)\\
\rho((a\ot b)\Delta(c))&=(a\ot b)\Delta(cm)
\end{align}
where $a,b,c\in A$. One can argue that these maps are well-defined, just as in the proof of Proposition \ref{prop:ext}. Next define $\lambda$ and $\rho$ on $A\ot A$ by $\lambda(x)=\lambda(Ex)$ and $\rho(x)=\rho(xE)$. If already $x\in E(A\ot A)$, the new definition coincides with the original one.
\ssnl
One can show again that this pair $(\lambda, \rho)$ is a multiplier. We denote it as $\Delta_1(m)$. Obviously $\Delta_1(1)=E$. It is also straightforward to show that $\Delta_1$ is a homomorphism extending $\Delta$. Finally, under the restriction $\Delta_1(1)=E$, this extension is unique.
\ebew

As before, we will again use $\Delta$, also for this extension $\Delta_1$.
\ssnl
Just as in the case of a non-degenerate homomorphism, this is a special case of the following more general result.

\prop
Let $A$ and $B$ be non-degenerate algebras and $\gamma:A\to M(B)$  a homomorphism. Assume that there is an idempotent $e\in M(B)$ with the property that $\gamma(A)B=eB$ and $B\gamma(A)=Be$. Then there is a unique homomorphism $\gamma_1:M(A)\to M(B)$ that extends $\gamma$ and satisfies $\gamma_1(1)=e$.
\eprop

The argument is found in \cite{VD-W1}. 
In fact, it is easier to give the proof in the more general situation.
\ssnl
The result can be applied to the homomorphisms $\Delta\ot\iota$ and $\iota\ot\Delta$ from $A\ot A$ to \\ 
$M(A\ot A\ot A)$, but only under the extra assumption that $A$ is idempotent. On the other hand, just as in the case of a regular coproduct, also here it is possible to adapt the proof along the same lines, to obtain the result without the assumption that $A$ is idempotent. 
\ssnl
We get the following.

\prop\label{prop:ext4}
Assume that $\Delta$ is a weakly non-degenerate homomorphism from $A$ to $M(A\ot A)$. Then $\Delta\ot\iota$ and $\iota\ot\Delta$ have unique extensions to homomorphisms from $M(A\ot A)$ to $M(A\ot A\ot A)$, denoted by the same symbols, satisfying $(\Delta\ot\iota)(1)=E\ot 1$ and $(\iota\ot\Delta)(1)=1\ot E$. 
\eprop

Here $1$ stands for the identity in $M(A)$ as well as for the identity in $M(A\ot A)$, while $E$ is the canonical idempotent as in Definition \ref{defin:wnd}.
\ssnl
This allows us again to give an alternative definition of coassociativity. 

\defin\label{defin:coass3}
Let $\Delta$ be a weakly non-degenerate homomorphism from $A$ to $M(A\ot A)$. We call $\Delta$ coassociative if $(\Delta\ot\iota)\Delta=(\iota\ot\Delta)\Delta$.
\edefin

Now $\Delta\ot\iota$ and $\iota\ot\Delta$ are the extensions obtained in Proposition \ref{prop:ext4}. The formula is required to hold on $A$ but just as in the case of a non-degenerate homomorphism, it holds on all of $M(A)$ as well, with the extension of $\Delta$ to $M(A)$ as obtained in Proposition \ref{prop:ext3}.
\ssnl
This notion is obviously compatible with the notion of coassociativity for a non-degenerate coproduct, as given in Definition \ref{defin:coass2} in the sense that, when $\Delta$ happens to be non-degenerate, i.e.\ when $E=1$, then the two notions coincide. This simply follows because we get the same extensions for $\Delta\ot\iota$ and $\iota\ot\Delta$. But again, we need to show that it is compatible with the original notions of coassociativity as in Definition \ref{defin:coass}.
More precisely, we need Proposition \ref{prop:equiv} also for weakly non-degenerate homomorphisms. However, this is again straightforward and the argument is precisely as in the proof of that proposition.

\opm
In this case, it does not seem to be true that weak non-degeneracy of the coproduct implies that $A$ is idempotent as soon as one of the canonical maps has range in $A\ot A$, a property that we could prove in the case of a non-degenerate product (in Proposition \ref{prop:idempotent}). To see this, assume that $T_1$ has range in $A\ot A$ and that $\Delta$ is weakly non-degenerate with canonical idempotent $E$. Assume that $\omega$ is a linear functional on $A$ that is $0$ on $A^2$. Then $(\omega\ot\iota)(\Delta(a)(b\ot c))=0$ for all $a,b,c\in A$. So $(\omega\ot\iota)(E(b\ot c))=0$ for all $b,c$. We can only conclude from this that $\omega$ is $0$ on all of $A$ if the left leg of $E$ is large enough. Indeed, clearly in the extreme case where $\Delta=0$ and so $E=0$, we can conclude nothing about $A$.
\eopm

Remark in passing that this is the reason why, in the theory of weak multiplier Hopf algebras, it is necessary to assume that the underlying algebra is idempotent, see Definition 1.14 in \cite{VD-W0}.

%
%

 \section{\hspace{-17pt}. Construction of the dual algebra} \label{s:dual} 
 
When $A$ is a vector space and $\Delta:A\to A\ot A$ a coassociative linear map, the linear dual $A'$ becomes an associative algebra if we define $\langle a,bb'\rangle=\langle \Delta(a),b\ot b'\rangle$ for $a\in A$ and $b,b'\in B$. We use the pairing notation for the evaluation of an element $b\in A'$ in a point $a\in A$, as well as for the tensor products. 
\ssnl
When $A$ is a non-degenerate algebra and $\Delta:A\to M(A\ot A)$ a coassociative linear map (as in the previous section), this is in general no longer possible as we can not apply $b\ot b'$ on elements $\Delta(a)$ for $b,b'\in A'$ and $a\in A$. Therefore, we can not make $A'$ into an algebra in this way. We need to consider suitable subspaces of $A'$ and also some natural conditions on $\Delta$.
\ssnl
In this section we will see what can be done. There are several possibilities.

\nl
\bf The algebra of reduced linear functionals \rm
\nl
In what follows, we assume that $A$ is a non-degenerate algebra and  $\Delta:A\to M(A\ot A)$  a linear map. We consider the canonical maps $T_1, T_2,T_3$ and $T_4$ 	as in Definition \ref{defin:1.2} of the previous section.

\defin\label{defin:2.1d}
Define $B^0_\ell$ as the subspace of linear functionals on $A$ spanned by elements of the form $f(c\,\cdot\,)$ where $f\in A'$ and $c\in A$.  Similarly we define $B^0_r$ as the space spanned by elements of the from $f(\,\cdot\,c)$ where $f\in A'$ and $c\in A$.
\edefin

Elements in $B^0_\ell$ and $B^0_r$ are sometimes called \emph{reduced linear functionals}. Remark that regularity of canonical maps is not needed to define these sets.
\ssnl
We can make these subspaces into  associative algebras under certain extra conditions on the coproduct. 
As we remarked before, there seems to be no way to define the product on the  full dual space $A'$ by the formula $(fg)(a)=(f\ot g)\Delta(a)$ for all $a$ because there is no way to define $f\ot g$ on $M(A\ot A)$ (in the general situation).
\ssnl
We now proceed in different steps below to make $B^0_\ell$ into an associative algebra.

\prop\label{prop:2.2d}
Assume that the canonical map $T_2$ is regular. If also $T_1$ or $T_3$ is regular,  we can define the product $\omega_1\omega_2$ for $\omega_1\in B^0_\ell$ and $\omega_2\in A'$ as a linear functional on $A$ so that  
\begin{equation}
(\omega_1\omega_2)(a)=\omega_2((f_1\ot \iota)((c_1\ot 1)\Delta(a)))\label{eqn:2.1c}
\end{equation}
when $\omega_1=f_1(c_1\,\cdot\,)$.
\eprop

\bew
We want to use the formula above to define this product. We just have to show that it is well-defined.
\ssnl
To do this, assume that $\sum_j f_j(c_j\,\cdot\,)=0$ where $f_j\in A'$ and $c_j\in A$. First assume that $T_1$ is regular. Then $\Delta(a)(1\ot c)$ belongs to $A\ot A$ for all $a,c$ and so
\begin{equation*}
\sum_j(f_j\ot\iota)((c_j\ot 1)\Delta(a)(1\ot c))=0
\end{equation*}
for all $a,c$. Because  $T_2$ is regular, we can cancel $c$ and obtain
\begin{equation*}
\sum_j(f_j\ot\iota)((c_j\ot 1)\Delta(a))=0
\end{equation*}
for all $a$. When we apply any linear functional $\omega_2$ of $A'$ we get 
$$\omega_2((f_1\ot \iota)((c_1\ot 1)\Delta(a)))=0.$$
This proves that the product is well-defined and that Equation (\ref{eqn:2.1c}) holds. 
\ssnl
In the case of regularity of $T_3$ we use that $(1 \ot c)\Delta(a)$ belongs to $A\ot A$ for all $a,c$ and hence
\begin{equation*}
\sum_j(f_j\ot\iota)((c_j\ot c)\Delta(a))=0
\end{equation*}
for all $a,c$. We can again cancel $c$ and proceed as in the previous case.
\ebew

We need regularity of $T_2$ for the Equation (\ref{eqn:2.1c}) to make sense. To prove that the product is well-defined, we used that  $T_1$ or $T_3$ is regular.

\opm
As a matter of fact, one could wonder if the regularity of  $T_2$ is sufficient to define $(\omega\ot\iota)\Delta(a)$ in $A$ when $\omega$ is of the form $f(c\,\cdot\,)$. However, it is not clear how to show in general that $(\omega\ot\iota)\Delta(a)$ is well-defined for $\omega$ in the space spanned by such functionals. 
\eopm

On the other hand, there is still another possibility to obtain that this product is well-defined as we see in the next proposition.

\prop
Assume that $A$ is an algebra with local units.
If $T_2$ is regular we can define the product $\omega_1\omega_2$ on $A$ when $\omega_1\in B^0_\ell$ and $\omega_2\in A'$ as in the previous proposition.
\eprop
\bew
To show that the product is well-defined, we start as in the proof of the previous proposition. 
So we  assume that 
$\sum_j f_j(c_j\,\cdot\,)=0$ where $f_j\in A'$ and $c_j\in A$. Now we use that $A$ has local units. Then we have an element $e\in A$ such that $c_i=c_ie$ for all $i$. We get, because $(e\ot 1)\Delta(a)\in A\ot A$ that
\begin{equation*}
\sum_j (f_j\ot\iota)((c_j\ot 1)\Delta(a))=
\sum_j (f_j(c_j\,\cdot)\ot\iota)((e\ot 1)\Delta(a))=0
\end{equation*}
and so also $\sum_j \omega_2((f_j\ot\iota)((c_j\ot 1)\Delta(a)))=0$ for any $\omega_2\in A'$. 
\ebew

We could continue with any of these possibilities. However, we restrict ourselves and we \emph{further assume that the canonical maps $T_1$ and $T_2$ are regular}. We proceed with the product as obtained in Proposition \ref{prop:2.2d} for this case. 

\prop\label{prop:2.5c}
If the coproduct $\Delta$ is a \emph{non-degenerate homomorphism} (as in Definition \ref{defin:nd}), then $\omega_1\omega_2\in B^0_\ell$ when $\omega_1,\omega_2$ are both in $B^0_\ell$. The result is still true if $\Delta$ is only weakly non-degenerate as in Definition \ref{defin:wnd}.
\eprop

\bew
i) Let $\omega_i=f_i(c_i\,\cdot\,)$ with $f_i\in A'$ and $c_i\in A$ for $i=1,2$. If we assume that  the coproduct $\Delta$ is a non-degenerate homomorphism as in Definition \ref{defin:nd}, we can write
\begin{equation*}
c_1\ot c_2=\sum_j (p_j\ot q_j)\Delta(r_j)
\end{equation*}
and then
\begin{align*}
(f_1\ot f_2)((c_1\ot c_2)\Delta(a))
&=\sum_j (f_1\ot f_2)((p_j\ot q_j)\Delta(r_j)\Delta(a))\\
&=\sum_j (f_1\ot f_2)((p_j\ot q_j)\Delta(r_ja)).
\end{align*}
This is equal to $\sum_j g_j(r_ja)$ when we define $g_j$ on $A$ by $g_j(x)= (f_1\ot f_2)((p_j\ot q_j)\Delta(x))$ for all $x\in A$. We see that $\omega_1\omega_2\in B^0_\ell$ when $\omega_i=f_i(c_i\,\cdot\,)$ for $i=1,2$.
\ssnl
ii) Now we show that the result is still true when $\Delta$ is weakly non-degenerate as in Definition \ref{defin:wnd}. Now we write 
\begin{equation*}
(c_1\ot c_2)E=\sum_j (p_j\ot q_j)\Delta(r_j)
\end{equation*}
where $E$ is the idempotent from the definition. We still get
\begin{equation*}
(c_1\ot c_2)\Delta(a)=\sum_j (p_j\ot q_j)\Delta(r_ja)
\end{equation*}
because $(c_1\ot c_2)\Delta(a)=(c_1\ot c_2)E\Delta(a)$. Then we can proceed as in item i).
\ebew

In the following proposition, we see when this product on $B^0_\ell$ is associative.

\prop\label{prop:2.4d}
Assume that $\Delta$ is a non-degenerate or a weakly non-degenerate homomorphism. Assume that it is coassociative in the sense of item i) of Definition \ref{defin:coass}. 
If $\omega_1, \omega_2\in B^0_\ell$ and $\omega_3\in A'$, then $(\omega_1\omega_2)\omega_3=\omega_1(\omega_2\omega_3)$. In particular, the product in $B^0_\ell$ is associative and we have a left action of the algebra $B^0_\ell$ on $A'$.
\eprop

\bew
i) Take $c,c_1,a$ in $A$. Because we assume that the maps $T_1$ and $T_2$ are regular and $\Delta$ coassociative, we find that
\begin{equation}
(c_1\ot 1\ot 1)((\Delta\ot\iota)(\Delta(a)(1\ot c)))=((\iota\ot\Delta)((c_1\ot 1)\Delta(a)))(1\ot 1\ot c).\label{eqn:2.2c}
\end{equation} 
Multiply with an element $c_2$ in the second factor from the left and then apply linear functionals $f_1$ and $f_2$ on the first and the second factor. For the left hand side we get, with $\omega_1=f_1(c_1\,\cdot\,)$ and $\omega_2=f_2(c_2\,\cdot\,)$,
\begin{align*}
(\omega_1\ot\omega_2\ot\iota)((\Delta\ot\iota)(\Delta(a)(1\ot c)))
&=(\omega_1\omega_2\ot\iota)(\Delta(a)(1\ot c))\\
&=((\omega_1\omega_2\ot\iota)\Delta(a))c.
\end{align*}
For the last step we use that $\omega_1\omega_2$ is in $B^0_\ell$ and so a linear combination of functionals of the form $f'(c'\,\cdot\,)$. For the right hand side of Equation (\ref{eqn:2.2c}) we obtain
\begin{align*}
(f_1\ot f_2\ot\iota)&((1\ot c_2\ot 1)(\iota\ot\Delta)((c_1\ot 1)\Delta(a))(1\ot 1\ot c))\\
&=((f_2(c_2\,\cdot\,)\ot \iota)\Delta((f_1(c_1\,\cdot\,)\ot\iota)\Delta(a)))c\\
&=((\omega_2\ot\iota)\Delta((\omega_1\ot\iota)\Delta(a)))c.
\end{align*}
We can cancel $c$ to get
\begin{equation*}
(\omega_1\omega_2\ot\iota)\Delta(a)=(\omega_2\ot\iota)\Delta((\omega_1\ot\iota)\Delta(a)).
\end{equation*}
\ssnl
ii) If we apply any $\omega_3$ of $A'$ we get $(\omega_1\omega_2)\omega_3=\omega_1(\omega_2\omega_3)$. With $\omega_3$ again in $B^0_\ell$ we get coassociativity of the product. Then, with any $\omega_3$ we obtain that $A'$ is a left $B^0_\ell$-module.
\ebew

Finally, we consider non-degeneracy of the product in $B^0_\ell$. It is clear that in a sense we need that the legs of $\Delta$ are all of $A$ in order to prove that the product induced from the coproduct is non-degenerate. Therefore, the following result is expected.

\prop Assume that the coproduct is full in the sense of Definition \ref{defin:1.16c}. Given $\omega_2\in A'$, then $\omega_2=0$ if $\omega_1\omega_2=0$ for all $\omega_1\in B_\ell$. On the other hand, given $\omega_1\in B_\ell$, then $\omega_1=0$ if $\omega_1\omega_2=0$ for all $\omega_2\in B_\ell$. In particular, the product in $B^0_\ell$ is non-degenerate.
\eprop

\bew 
i) Assume that $\omega_2\in A'$ and that $\omega_1\omega_2=0$ for all $\omega_1\in B^0_\ell$. If $\omega_1=f(c\,\cdot\,)$ for $f\in A'$ and $c\in A$ we get by the definition of the product that
\begin{equation*}
0=(\omega_1\omega_2)(a)=\omega_2((f(c\,\cdot\,)\ot\iota)\Delta(a)).
\end{equation*}
This holds for all $a,c\in A$ and $f\in A'$. Because $\Delta$ is assumed to be full, it follows that $\omega_2(x)=0$ for all $x\in A$ and so $\omega_2=0$.
\ssnl
ii) Now assume that $\omega_1$ is in $B^0_\ell$ and that $\omega_1\omega_2=0$ for all $\omega_2\in B^0_\ell$. Then

$
f(c((\omega_1\ot \iota)\Delta(a)))=0
$ 
for all $f$ and $c$. Then also $c((\omega_1\ot \iota)\Delta(a))=0$ for all $c$. We can cancel $c$, multiply with $c'$ on the other side to obtain $(\omega_1\ot \iota)\Delta(a))c'=0$. Apply again any $f$ and we get $\omega_1((\iota\ot f(\, \cdot\,c'))\Delta(a))=0$ for all $a,c'$ and $f$. Again because $\Delta$ is full we get $\omega_1(x)=0$ for all $x$ and so $\omega_1=0$.
\ebew

\oldcomment{What about this:
Assume that  $f\in A'_r$ and that $fg=0$ with $g=\omega(\,\cdot\,c)$ for all $\omega$ and all $c$. This means that $(f\ot\iota)(\Delta(a)(1\ot c))=0$ for all $a,c$. Because we have that  $T_1$ is  surjective, it follows that $f=0$. On the other hand, take again $f\in A'_r$ and now assume that $gf=0$ with $g=\omega(\,\cdot\,c)$ for all $\omega$ and all $c$. This will imply that $(\iota\ot f)(\Delta(a)(c\ot 1))=0$ for all $a$ and all $c$. In order to conclude from this that again $f=0$, we need some other assumptions. A possibility is that also $T_4$ is assumed to be regular and surjective. We can cancel $c$ on the right and multiply again on the left. Then we see that another possibility is that $T_2$ is regular and surjective. Finally, also the existence of a counit would do the job.
}{}

We summarize the results obtained so far. \oldcomment{Check it again!}

\prop\label{stel:2.6}
Let $A$ be a non-degenerate algebra and $\Delta: A\to M(A\ot A)$ a non-degenerate or weakly non-degenerate homomorphism. Assume that  the canonical maps $T_1$ and $T_2$ are regular and that $\Delta$ is a full coproduct. Finally require that it is coassociative in the sense of Definition \ref{defin:coass}. Then the space $B^0_\ell$ of linear functionals on $A$, spanned by elements of the form $f(c\,\cdot\,)$ where $c\in A$ and $f\in A'$ is a non-degenerate associative algebra for the product defined by 
\begin{equation*}
(\omega_1\omega_2)(a)=\omega_2((\omega_1\ot\iota)\Delta(a)). 
\end{equation*}
The space $A'$ is a left $B^0_\ell$-module for the action $(\omega_1,\omega_2)\mapsto  \omega_1\omega_2$ where we use the same formula as above.
\eprop
We have a similar result for $B^0_r$, defined  as the space spanned by elements of the from $f(\,\cdot\,c)$ where $f\in A'$ and $c\in C$:

\prop \label{stel:2.7}
With the same conditions as in the previous theorem, the space $B^0_r$ of linear functionals on $A$, spanned by elements of the form $f(\,\cdot\,c)$ where $c\in A$ and $f\in A'$ is a non-degenerate associative algebra for the product defined by 
\begin{equation*}
(\omega_1\omega_2)(a)=\omega_1((\iota\ot \omega_2)\Delta(a)). 
\end{equation*}
The space $A'$ is a right $B^0_r$-module for the action $(\omega_1,\omega_2)\mapsto  \omega_1\omega_2$ where we use the same formula as above.
\eprop

\oldcomment{Do we not have that local units exist under these conditions? Check and make a remark about the possibility of the other approach.}{}

We can consider the intersection $B^0_\ell\cap B^0_r$. We denote it by $B_0$. For two elements $\omega_1$ and $\omega_2$ we have the product on $B^0_\ell$ as in Proposition  \ref{stel:2.6} and the one on $B^0_r$ as in Proposition \ref{stel:2.7}

\prop \label{prop:2.8}
If $\omega_1$ and $\omega_2$ belong to $B_0$, then the two products coincide. Therefore $B_0$ is a subalgebra of both $B^0_\ell$ and $B^0_r$. 
\eprop
\bew
i) For the product $\omega_1\omega_2$ in $B^0_\ell$ we have, for $a\in A$, 
\begin{equation*}
(\omega_1\omega_2)(a)=\omega_2((\omega_1\ot\iota)\Delta(a)).
\end{equation*}
Take $\omega_1=f_1(c_1\,\cdot\,)$ and $\omega_2=f_2(\,\cdot\,c_2)$ where $f_1,f_2\in A'$ and $c_1,c_2\in A$. Then we have
\begin{align*}
\omega_2((\omega_1\ot\iota)\Delta(a))
&=f_2(((\omega_1\ot\iota)\Delta(a))c_2)\\
&=f_2((\omega_1\ot\iota)(\Delta(a)(1\ot c_2)))\\
&=f_2((f_1\ot\iota)((c_1\ot 1)\Delta(a)(1\ot c_2)))\\
&=(f_1\ot f_2)((c_1\ot 1)\Delta(a)(1\ot c_2))
\end{align*}
ii) In a similar way, we get for the product in $B^0_r$, for all $a\in A$, also
\begin{equation*}
\omega_1((\iota\ot\omega_1)\Delta(a))
=(f_1\ot f_2)((c_1\ot 1)\Delta(a)(1\ot c_2))
\end{equation*}
if $\omega_1=f_1(c_1\,\cdot\,)$ and $\omega_2=f_2(\,\cdot\,c_2)$.
\ssnl
iii) It is a straightforward consequence that $B$ is a subalgebra of the algebras $B^0_\ell$ and $B^0_r$.
\ebew

We have seen that the products in $B^0_\ell$ and $B^0_r$ are non-degenerate. It is also true for the subalgebra $B_0$ as we see in the next proposition.

\prop\label{prop:2.11d} 
If the coproduct is full, the product in $B_0$ is non-degenerate.
\eprop
\bew
i) Assume that $\omega_1\in B_0$ and that $\omega_1\omega_2=0$ for all $\omega_2\in B_0$. This means that $\omega_2((\omega_1\ot\iota)\Delta(a))=0$ because $B_0\subseteq B_\ell$. We can take for $\omega_2$ elements of the form $f(c\,\cdot\,c')$ with $f\in A'$ and $c,c'\in A$. It follows that $(\omega_1\ot\iota)\Delta(a)=0$ for all $a$. By the fullness of the coproduct, this implies that $\omega_1=0$.
\ssnl
ii) On the other hand assume that $\omega_2\in B_0$ and that $\omega_1\omega_2=0$ for all $\omega_1\in B_0$. This means that $\omega_1((\iota\ot\omega_2)\Delta(a))=0$ because $B_0\subseteq B_r$. As above, this implies that $(\iota\ot\omega_2)\Delta(a)=0$
for all $a$ and so also now $\omega_2=0$.
\ebew

Remark that the span of elements of the form $f(c\,\cdot\,c')$ with $f\in A'$ and $c,c'\in A$ give a subalgebra of $B_0$. This subalgebra is considered in the original paper on multiplier Hopf algebras, see Section 6 in \cite{VD-mha}. An even smaller subalgebra is considered in the case of a multiplier Hopf algebra with integrals, see \cite{VD-alg}. See also Remark \ref{opm:2.12d}.
This is of a different nature and beyond the scope of this paper.
\nl
If the coproduct admits a counit as in Definition \ref{defin:1.15c} one can expect that these algebras have a unit. We have the following result.

\prop\label{prop:2.11c}
If there is a counit $\varepsilon$ that is a homomorphism, the algebras $B^0_\ell$ and $B^0_r$, and consequently also $B_0$, are unital.
\eprop
\bew
i) If $c$ is an element of $A$ satisfying $\varepsilon(c)=1$ we see that 
\begin{equation*}
\varepsilon(ca)=\varepsilon(c)\varepsilon(a)=\varepsilon(a)
\end{equation*}
So $\varepsilon=\varepsilon(c\,\cdot\,)$ and $\varepsilon\in B^0_\ell$. Similarly $\varepsilon\in B^0_r$ and therefore also $\varepsilon\in B_0$.
\ssnl
ii) Let $\omega$ be an element in $B^0_\ell$ of the form $f(c\,\cdot\,)$. Then, using the product as defined in Proposition \ref{prop:2.2d}, 
\begin{equation*}
(\omega\varepsilon)(a)=(f\ot\varepsilon)((c\ot 1)\Delta(a))=f(ca)=\omega(a).
\end{equation*}
On the other hand
\begin{equation*}
(\varepsilon\omega)(a)=(\varepsilon\ot\omega)\Delta(a)=\omega(a)
\end{equation*}
because also $(\varepsilon\ot\iota)\Delta(a)=a$. This last property can be seen by applying $\varepsilon\ot \iota$ on $(c\ot 1)\Delta(a)$ with $\varepsilon(c)=1$. We see that $\varepsilon$ is a unit in the algebra $B^0_\ell$.
\ssnl
iii) In a similar way, $\varepsilon$ is a unit in the algebra $B^0_\ell$ and as it belong to $B_0$, it is also a unit for $B_0$.
\ebew

It is not clear if it is necessary to have that the counit is a homomorphism. We come back to this in the next item, see Remark 
\ref{opm:2.17c}. 

\opm\label{opm:2.12d}
Before we continue, let us briefly comment on the special cases where there is an integral (like for algebraic quantum groups, algebraic quantum hypergroups and algebraic quantum groupoids). In these cases, one takes for the dual algebra a space of functionals of the form $\varphi(\,\cdot\, c)$ where $c\in A$ and $\varphi$ an integral. One uses the properties of integrals to show that we get an algebra of linear functionals. We will not consider these cases further.
We refer to the literature on these objects for results of this type. 
\eopm

In a further item, we will discuss the pairing of $A$ with these algebras $B^0_\ell$, $B^0_r$ and $B_0$.
\ssnl
First, we look for a bigger subspace of $A'$ that still can be made into an associative algebra under certain natural conditions on the coproduct.

\oldcomment{\ssnl Referenties naar de litteratuur? \rood To do later.}{}

\nl
\bf The dual algebra - A different approach \rm
\nl
Next we will try a different approach. 
It has some advantages over the previous one, but also some disadvantages as we will see. 
\ssnl
As before, we have a non-degenerate algebra $A$ and a linear map $\Delta:A\to M(A\ot A)$ and we consider the canonical maps $T_1$, $T_2$, $T_3$ and $T_4$.

\defin\label{defin:2.6}
i) Assume that the canonical map $T_1$ is regular. Then for all $a\in A$ and $\omega\in A'$ we can define a left multiplier $x$ of $A$ by
\begin{equation*}
xc=(\omega\ot\iota)(\Delta(a)(1\ot c))
\end{equation*}
for $c\in A$.  We denote this multiplier by $(\omega\ot \iota)\Delta(a)$. We  define $B_\ell$ as the set of elements $\omega\in A'$ such that $(\omega\ot\iota)\Delta(a)$ is an element of $A$ for all $a\in A$.
\ssnl
ii) Similarly, when $T_2$ is regular, we have a right multiplier $(\iota\ot\omega)\Delta(a)$ for all $a\in A$ and $\omega\in A'$ given by
\begin{equation*}
cx=(\iota\ot\omega)((c\ot 1)\Delta(a))
\end{equation*}
for $c\in A$.  We  denote by $B_r$ the space of elements $\omega\in A'$ such that this multiplier belongs to $A$. 
\edefin

We could also consider the other canonical maps $T_3$ and $T_4$ to define these sets. In the event that $T_1$ and $T_3$ are regular, the two notions  coincide.  In fact then $(\omega\ot\iota)\Delta(a)$ will be an element in $M(A)$ for all $a$ and if $\omega\in B_\ell$, this will be in $A$ for all $a$. Similarly if $T_2$ and $T_4$ are regular, $(\iota\ot\omega)\Delta(a)\in M(A)$ for all $a$ and if $\omega\in B_r$, it is an element of $A$ for all $a$.
\ssnl
 In general, we have no information about the size of these spaces $B_r$ and $B_\ell$ here. If e.g.\ $\Delta$ is defined as $\Delta(a)=a\ot 1$, then $T_1$ is regular, but $B_\ell$ is trivial if $A$ has no unit. We also have no immediate relation between $B_\ell$ and $B_r$ in the event that both $T_1$ and $T_2$ are regular. 
\ssnl
As in the previous item, again in what follows  \emph{we will stick to the case} with $T_1$ and $T_2$ regular as in Definition \ref{defin:2.6} above.  But it is always good to have in mind that there are other possibilities. There are examples of this type in the next section.
\ssnl
We have the following set of inclusions.

 \prop\label{prop:2.15e}
 i) Assume that  $T_1$ is regular. If $T_2$ or $T_4$ is regular, then $B^0_r\subseteq B_r$.\\
 ii) Asume that $T_2$ is regular. If $T_1$ or $T_3$ is regular, then $B^0_\ell\subseteq B_\ell$.\\
 iii) Assume that $T_3$ is regular. If $T_2$ or $T_4$ is regular, then $B^0_\ell\subseteq B_r$. \\
 iv) Assume that $T_4$ is regular. If $T_1$ or $T_3$ is regular, then $B^0_r\subseteq B_\ell$
\eprop
\bew
Assume that $T_1$ is regular and that $\omega=f(\,\cdot\,c)$ where $f\in A'$ and $c\in A$. If $T_2$ is regular we have for all $c'\in A$,
\begin{equation*}
(\iota\ot\omega)((c'\ot 1)\Delta(a))=(\iota\ot f)((c'\ot 1)\Delta(a)(1\ot c))=c'x
\end{equation*}
where $x=(\iota\ot f)(\Delta(a)(1\ot c))$. Therefore $\omega\in B_r$. Similarly if $T_4$ is regular we have for all $c'\in A$, 
\begin{equation*}
(\iota\ot\omega))(\Delta(a)(c'\ot 1))=( \iota\ot f)(\Delta(a)(c'\ot c))=xc'
\end{equation*}
where $x=(\iota\ot f)(\Delta(a)(1\ot c))$. Therefore again $\omega\in B_r$.
\ssnl
The other cases are proven in the same way.
\ebew

In the case where $T_1$ and $T_2$ are regular,  the spaces $B_\ell$ and $B_r$  we define in Definition \ref{defin:2.6}  contain the spaces defined earlier in Definition  \ref{defin:2.1d}, Proposition  \ref{stel:2.6} and Proposition \ref{stel:2.7}.
\nl
We now look for conditions that allow us to make also these bigger  spaces $B_\ell$ and $B_r$  into associative algebras. We proceed as for the more restricitive spaces in the previous item.
\ssnl
A first step is to define the product in $B_\ell$, dual to the coproduct on $A$. We generalize here the result of Proposition \ref{prop:2.2d}

\defin\label{defin:2.7a}
i) We define a bilinear map $(\omega_1,\omega_2)\mapsto \omega_1\omega_2 $ from $B_\ell \times A'$ to $A'$ by
\begin{equation*}
(\omega_1\omega_2)(a)=\omega_2((\omega_1\ot\iota)\Delta(a))
\end{equation*}
where $a\in A$. This is possible because  $(\omega_1\ot\iota)\Delta(a)\in A$ for $\omega_1\in B_\ell$ by definition. 
\ssnl
ii) Similarly, we have a bilinear map $(\omega_1,\omega_2)\mapsto \omega_1\omega_2 $ form $A'\times B_r$ to $A'$ given by
\begin{equation*}
(\omega_1\omega_2)(a)=\omega_1((\iota\ot\omega_2)\Delta(a))
\end{equation*}
where $a\in A$. This is possible because $(\iota\ot\omega_2)\Delta(a)$ is well-defined in $A$ when $\omega_2\in B_r$.
\edefin

We should in fact be more careful with the notations here. If $\omega_1\in B_\ell$ and $\omega_2\in B_r$ we have
\begin{align*}
(\omega_1\omega_2)(a)&=\omega_2((\omega_1\ot\iota)\Delta(a)) \\
(\omega_1\omega_2)(a)&=\omega_1((\iota\ot\omega_2)\Delta(a))
\end{align*}
for all $a$. For the first formula, we use that $\omega_1\in B_\ell$ whereas for the second one that $\omega\in B_r$. It is not clear whether or not these two give the same result. We will show later that this is the case under certain conditions, see Proposition \ref{prop:2.10b}. \oldcomment{Add a reference!}{}
\ssnl
This remark is only important if we consider the two products at the same time. For the moment, this is not  the case and so the problem does not yet occur.

\opm\label{opm:2.17c}
i) If there is a counit $\varepsilon$ we have $(\varepsilon\ot\iota)\Delta(a)=a$ and so $\varepsilon\omega=\omega$ for all $\omega\in A'$. We are using here the product on $B_\ell\times A'$. Similarly we have $\omega\varepsilon=\omega$ for the product on $A'\times B_r$. 
\ssnl
ii) It is not clear if also $\omega\varepsilon=\omega$ when $\omega\in B_\ell$ and $\varepsilon\omega=\omega$ for $\omega\in B_r$. To have the first property, we would need that $\varepsilon((\omega\ot\iota)\Delta(a))=\omega(a)$ for $\omega\in B_\ell$. This is the case when $\varepsilon$ is a homomorphism. Indeed, if  $\varepsilon(c)=1$  then
\begin{align*}
\varepsilon((\omega\ot\iota)\Delta(a))
&=\varepsilon(((\omega\ot\iota)\Delta(a))c) \\
& =\varepsilon((\omega\ot\iota)\Delta(a)(1\ot c)) \\
&=\omega((\iota\ot\varepsilon)(\Delta(a)(1\ot c)))=\omega(a).
\end{align*}

\eopm

\oldcomment{Refer to a later property}{\ssnl}

\nl
\bf A stronger form of coassociativity \rm
\nl
We will need a new form of coassociativity. To introduce this, assume that the maps $T_1$ and $T_2$ are regular and that $\Delta$ is coassociative as in item i) of Definition \ref{defin:coass}. So, for all $a,c,c'$ in $A$ we have 
\begin{equation}
(c\ot 1\ot 1)((\Delta\ot\iota))(\Delta(a)(1\ot c'))
=((\iota\ot\Delta)((c\ot 1)\Delta(a))(1\ot 1\ot c').\label{eqn:2.1a}
\end{equation}
We proceed as in the proof of Proposition \ref{prop:2.4d}. 
Let $f,g\in A'$ and denote $f(c\,\cdot\,)$ by $\omega_1$ and $g(\,\cdot\,c')$ by $\omega_2$. If we apply $f\ot\iota\ot g$ on Equation (\ref{eqn:2.1a}) we find for all $a\in A$,
\begin{equation*}
(\omega_1\ot\iota)(\Delta((\iota\ot\omega_2)\Delta(a)))
=(\iota\ot\omega_2)(\Delta((\omega_1\ot\iota)\Delta(a))).
\end{equation*}
This leads to the following \emph{new notion of coassociativity} involving the spaces $B_\ell$ and $B_r$. 

\defin\label{defin:2.8}
Assume that the maps $T_1$ and $T_2$ are regular. Then we call $\Delta$ coassociative if
\begin{equation}
(\omega_1\ot\iota)(\Delta((\iota\ot\omega_2)\Delta(a)))
=(\iota\ot\omega_2)(\Delta((\omega_1\ot\iota)\Delta(a)))\label{eqn:2.2b}
\end{equation}
for all $a\in A$ and all $\omega_1\in B_\ell$ and $\omega_2\in B_r$.
\edefin

We see from the previous observation that this notion is (at least in principle) stronger than the one formulated in  Definition \ref{defin:coass}. 
\ssnl
We get a similar property for the case where $T_3$ and $T_4$ are regular. Then we consider coassociativity as in item ii) of Definition \ref{defin:coass}.

\opm
Once more, we see that we are using various forms of coassociativity of the coproduct. This may be confusing. However, instead of giving these different notions all a name, we systematically speak about \emph{coassociativity in the sense of} and refer to the formulation of the notion we are actually using (as we have done before).
\eopm
\keepcomment{We must verify that we do this! Also in the examples. \rood To do!}{}

One of the consequences of the above stronger form of coassociativity is the following result.

\oldcomment{\ssnl 
Refer to this earlier and also refer to this earlier statement here.}{}

\prop\label{prop:2.10b}
Suppose that $T_1$ and $T_2$ are regular and that the coproduct is coassociative as in Definition \ref{defin:2.8} above. Also suppose that we have a counit $\varepsilon$ and that it is a homomorphism. Then 
\begin{equation*}
\omega_1((\iota\ot\omega_2)\Delta(a))=\omega_2((\omega_1\ot\iota)\Delta(a))
\end{equation*}
for all $\omega_1\in B_\ell$ and $\omega_2\in B_r$. 
\eprop

\bew
i) For $\omega_1\in B_\ell$ we have 
\begin{equation*}
\varepsilon((\omega_1\ot\iota)(\Delta(a)(1\ot c)))=\omega_1((\iota\ot\varepsilon)(\Delta(a)(1\ot c)))=\omega_1(a)\varepsilon(c)
\end{equation*}
but also
\begin{equation*}
\varepsilon((\omega_1\ot\iota)(\Delta(a)(1\ot c))=\varepsilon((\omega_1\ot\iota)\Delta(a))\varepsilon(c)
\end{equation*}
so that $\varepsilon((\omega_1\ot\iota)\Delta(a))=\omega_1(a)$.
\ssnl
ii) Similarly $\varepsilon((\iota\ot\omega_2)\Delta(a)=\omega_2(a)$ when $\omega_2\in B_r$.
\ssnl
iii) Therefore, if we apply $\varepsilon$ on Equation (\ref{eqn:2.2b}) we get the desired result.
\ebew

\keepcomment{
Can we prove it in another way (also when the counit is not an homomorphism, but when $\Delta$ is full? This notion may be usefull later also. \rood To do. \blauw Refer to this problem in the conclusions section and the examples.}{}

Suppose that we have given this result. If $\varepsilon$ is a counit, it would follow from the fact that $\varepsilon\in B_\ell$  that
\begin{equation*}
\varepsilon((\iota\ot\omega_2)\Delta(a))=\omega_2((\varepsilon\ot\iota)\Delta(a))=\omega_2(a)
\end{equation*}
for all $\omega_2\in B_r$. Similarly we would obtain that $\varepsilon((\omega_1\ot\iota)\Delta(a))=\omega_1(a)$ for $\omega_1\in B_r$. We see from the proof that these properties are not sufficient to get the result. We need that the counit is a homomorphism.

\oldcomment{ \ssnl We have to check this aspect on the whole. This is more or less done.}{\ssnl}

Before we continue, we make one more comment on this form of coassociativity.

\opm
i) If the canonical maps $T_2$ and $T_3$ are regular, we can define both \\ $(\omega\ot\iota)\Delta(a)$  and $(\iota\ot\omega)\Delta(a)$ as right multipliers. Then we can consider the space $B_\ell$ as the set of functionals $\omega$ such that  $(\omega\ot\iota)\Delta(a)\in A$ for all $a$ and $B_r$ as the set of functionals $\omega$ such that $(\iota\ot\omega)\Delta(a)$ belongs to $A$. Assume now that $\Delta$ is coassociative in the sense of item ii) of Definition \ref{defin:1.11a} . 
We can apply $f\ot 1\ot g$ on 
\begin{equation*}
(c\ot 1\ot 1)(\Delta\ot\iota)((1\ot c')\Delta(a))
=(1\ot 1\ot c')(\iota\ot\Delta)((c\ot 1)\Delta(a)).
\end{equation*}
If now $\omega_1=f(c\,\cdot\,)$ and $\omega_2=g(c'\,\cdot\,)$ we again get
\begin{equation*}
(\omega_1\ot\iota)(\Delta((\iota\ot\omega_2)\Delta(a)))
=(\iota\ot\omega_2)(\Delta((\omega_1\ot\iota)\Delta(a)))
\end{equation*}
for all $a\in A$ and all $\omega_1\in B_\ell$ and $\omega_2\in B_r$.
\ssnl
ii) We have a similar observation when $T_1$ and $T_4$ are regular. 
\ssnl
In all these cases, we get a similar expression for coassociativity.
\eopm
\keepcomment{\rood We have to add more remarks about this. ?? \blauw Refer to a statement in the previous section. And to statements later in this section. One is where we try to avoid this stronger form of coassociativity. See \ref{??}.
}{}

Now we show that the spaces $B_\ell$ and $B_r$, defined in  Definition \ref{defin:2.6}, are associative algebras. The proof is very similar as the one of Proposition \ref{prop:2.4d}.
\oldcomment{We should give the proof for $B_\ell$.}

\prop\label{prop:2.16b}
Assume that the maps $T_1$ and $T_2$ are regular and that $\Delta$ is coassociative in the sense of Definition \ref{defin:2.8}.
For $\omega_1, \omega_2$ both in $B_\ell$ we have $\omega_1\omega_2\in B_\ell$ making $B_\ell$ into an associative algebra. Moreover, $A'$ is a left $B_\ell$-module.
\eprop

\bew
i) Let $\omega_1\in B_\ell$ and $\omega_2\in A'$. Recall that we have defined the product $\omega_1\omega_2$ in $A'$ by $(\omega_1\omega_2)(a)=\omega_2((\omega_1\ot \iota)\Delta(a))$ for all $a\in A$, see Definition \ref{defin:2.7a} . 
Now take $f\in A'$ and $c\in A$ and denote $f(\,\cdot\,c)$ by $\omega_3$. Then we have for all $a\in A$,
\begin{align*}
f((\omega_1\omega_2\ot\iota)(\Delta(a)(1\ot c)))
&=(\omega_1\omega_2)((\iota\ot\omega_3)\Delta(a))\\
&=\omega_2(((\omega_1\ot\iota)\Delta)((\iota\ot\omega_3)\Delta(a))).
\end{align*}
By coassociativity we have
\begin{equation*}
((\omega_1\ot\iota)\Delta)((\iota\ot\omega_3)\Delta(a))
=((\iota\ot\omega_3)\Delta)((\omega_1\ot\iota)\Delta(a)).
\end{equation*}
Then we find
\begin{align*}
f((\omega_1\omega_2\ot\iota)(\Delta(a)(1\ot c)))
&=\omega_2(((\iota\ot\omega_3)\Delta)((\omega_1\ot\iota)\Delta(a)))\\
&=\omega_2((\iota\ot f)(\Delta((\omega_1\ot\iota)\Delta(a))(1\ot c))).
\end{align*}
If also $\omega_2\in B_\ell$ we can write
\begin{equation*}
f((\omega_1\omega_2\ot\iota)(\Delta(a)(1\ot c)))=f(((\omega_2\ot\iota)\Delta((\omega_1\ot \iota)\Delta(a)))c).
\end{equation*}
The equation holds for all $f$ and all $c$ and so we obtain that $(\omega_1\omega_2\ot\iota)\Delta(a)\in A$ and 
\begin{equation}
((\omega_1\omega_2)\ot\iota)\Delta(a)=(\omega_2\ot\iota)\Delta((\omega_1\ot \iota)\Delta(a)).\label{eqn:2.4d}
\end{equation}
This proves that $\omega_1\omega_2\in B_\ell$. 
\ssnl
ii) If we apply any other element $\omega$ of $B_\ell$ we obtain
\begin{equation*}
((\omega_1\omega_2)\omega)(a)=(\omega_2\omega)((\omega_1\ot\iota)\Delta(a))=(\omega_1(\omega_2\omega))(a)
\end{equation*}
and we see that the product in $B_\ell$ is associative. 
\ssnl
iii) If we apply any element $\omega$ of $A'$  to Equation (\ref{eqn:2.4d}) we get also
 \begin{equation*}
((\omega_1\omega_2)\omega)(a)=(\omega_1(\omega_2\omega))(a)
\end{equation*}
And we see that $A'$ is a left $B_\ell$-module. 
\ebew

\oldcomment{We are in fact using a weaker form of coassociativity to get Equation (\ref{eqn:2.4d}) but still stronger than the one in item ii) of Definition \ref{defin:1.11a}. On the other hand, we need the stronger version for this last property.}{}
\oldcomment{
\ssnl We have a form of coassociativity that is in fact just associativity of the product in the dual. \rood Clarify!}{}
\ssnl
We have a similar result for $B_r$.

\prop\label{prop:2.17b}
Assume that the maps $T_1$ and $T_2$ are regular and that $\Delta$ is coassociative in the sense of Definition \ref{defin:2.8}.
For $\omega_1, \omega_2$ both in $B_r$ we have $\omega_1\omega_2\in B_r$ making $B_r$ into an associative algebra. Moreover, $A'$ is a right $B_r$-module.
\eprop
\bew
The proof is completely similar as the one of the previous result. Now we work with functionals $f(c\,\cdot\,)$ to show that 
$\omega_1\omega_2\in B_r$ when $\omega_1,\omega_2\in B_r$ and that 
\begin{equation*}
(\iota\ot(\omega_1\omega_2)\Delta(a)=(\iota\ot\omega_1)\Delta((\iota\ot\omega_2)\Delta(a)).
\end{equation*}
Then we continue as in step ii) and iii) of the proof of the previous result.
\ebew
\keepcomment{We need to say something about coassociativity here and how we use it here. It is under a different form. See an earlier remark. This has been done.}{\ssnl}

Again, in order to have that this product is non-degenerate, we need that  $\Delta$ is full.

\oldcomment{
\ssnl We need to recall the notion of a full coproduct, with references. This is done. But we need to adapt this argument. Refer to the previous case}{}

\prop\label{prop:2.24d}
Assume that the canonical maps $T_1$ and $T_2$ are regular. If the coproduct is full, then the product in $B_\ell$ is non-degenerate and the action of $B_\ell$ on $A'$ is faithful.
\eprop
\bew
i) First assume that $f\in A'$ and that $\omega f=0$ for all $\omega\in B_\ell$. This means that $f((\omega\ot\iota)\Delta(a))=0$. We can take $\omega$ of the form $g(c\,\cdot\,)$. By the fullness of the coproduct all elements in $A$ are of the form $(\omega\ot\iota)\Delta(a)$ with such elements $\omega$ (see a remark after Definition \ref{defin:1.16c}). Therefore $f=0$.
\ssnl
ii) Next assume that $\omega_1\omega_2=0$ for all $\omega_2\in B_\ell$. Because there are enough elements in $B_\ell$ in the sense that $a=0$ if  $\omega(a)=0$ for all $\omega\in B_\ell$ we have that $(\omega_1\ot\iota)\Delta(a)=0$ for all $a$. Now multiply with any element $c$ from the right and apply any $f$ from $A'$. Because $\Delta$ is assumed to be full, again elements of the form
$(\iota\ot f(\,\cdot\,c))\Delta(a)$ span all of $A$. Therefore $\omega_1=0$.
\ebew

We have a similar result for $B_r$.
\ssnl
\keepcomment{This is also usable in the case of weak multiplier Hopf algebras. Add a reference}{}
\keepcomment{\ssnl What about the actions? We keep this for later. See file artikel2ra.tex}{\ssnl}

As we did in the previous item, with the algebras of reduced functionals (see Proposition \ref{prop:2.8}), we now also consider the intersection $B$, here defined as $B_r\cap B_\ell$. We get the following result.

\prop\label{stel:2.25c}
Let $A$ be a non-degenerate algebra and $\Delta$ a linear map from $A$ to $M(A\ot A)$. Assume that $T_1$ and $T_2$ are regular and that $\Delta$ is coassociative in the sense of Definition \ref{defin:2.8}.  Consider the algebras $B_r$ and $B_\ell$ as in Propositions \ref{prop:2.16b} and \ref{prop:2.17b}. If $\omega_1((\iota\ot\omega_2)\Delta(a))=\omega_2((\omega_1\ot\iota)\Delta(a)$ for all $\omega_1\in B_\ell$ and $\omega_2\in B_r$ and all $a\in A$, then the products on $B$ coincide. In particular, $B$ is a subalgebra of the algebra $B_\ell$ and of the algebra $B_r$. If the coproduct is full, the product in $B$ is still non-degenerate.
\eprop
\bew
i) Let $\omega_1,\omega_2\in B$. Because $\omega_1\in B_\ell$ we can use the product $\omega_1\omega_2$ as defined in item i) of Definition \ref{defin:2.7a}, given by $(\omega_1\omega_2)(a)=\omega_2((\omega_1\ot\iota)\Delta(a))$. Because $\omega_2\in B_r$ we can use the product $\omega_1\omega_2$ as defined in item ii) of Definition \ref{defin:2.7a}; given by $(\omega_1\omega_2)(a)=\omega_1((\iota\ot\omega_1)\Delta(a))$. By assumption these two expressions coincide and define $(\omega_1\omega_2)(a)$.
\ssnl
ii) Because $B_\ell$ is an algebra for this product (see Proposition \ref{prop:2.16b} and also $B_r$ is an algebra (see Proposition \ref{prop:2.17b}, we have $\omega_1\omega_2\in B$.
\ssnl
iii) The non-degeneracy of the product in $B$ here is a consequence of the fact that the pairing of $A$ with $B$ is non-degenerate. Indeed, suppose that $\omega_1\in B$ and that $\omega_1\omega_2=0$ for all $\omega_2\in B$. This means that $\omega_2((\omega_1\ot\iota)\Delta(a))=0$. Because there are enough elements in $B$, we get $(\omega_1\ot \iota)\Delta(a)=0$.  Then we can proceed as before to get from the fullness of the coproduct that $\omega_1=0$. Similarly when $\omega_2\in B$ and  $\omega_1\omega_2=0$ for all $\omega_1\in B$. In this case, we use the formula $\omega_1((\iota\ot \omega_2)\Delta(a))$ for $\omega_1\omega_2(a)$. 
\ebew
Remark that the proof is essentially the same as the one used to show that the algebra $B_0$ is non-degenerate in Proposition \ref{prop:2.11d}.

\prop
If there is a counit $\varepsilon$, it is a unit in the algebra $B$.
\eprop

\bew
If there is a counit $\varepsilon$, it belongs to the spaces $B_\ell$ and $B_r$ (see Remark \ref{opm:2.17c}. So it belongs to $B$. We also have seen in Remark \ref{opm:2.17c} that $\varepsilon\omega=\omega$ and  $\omega\varepsilon=\omega$  for all $\omega\in A'$.  In particular, this holds for all $\omega\in B$.
\ebew

Under the given conditions, namely (1) the more general form of coassociativity (as in Definition \ref{defin:2.8}) and (2) the equality of the two products, we get an algebra $B$ that is in general strictly bigger than the algebra we get in Proposition \ref{prop:2.8} of the previous item. On the other hand,  remember that in order to get the smaller algebras, although we only need a weaker form of coassociativity and the equality of the two products is obvious,  it is needed that $\Delta$ is a homomorphism, non-degenerate or weakly non-degenerate.
\ssnl
In many of the examples we consider in Section \ref{s:examp} we will see that the conditions are satisfied and hence we get this bigger algebra.

\keepcomment{Formulate as a result, check the argument and make a comment about the condition that $\varepsilon$ is a homomorphism. Compare with the previous result. Refer to the examples.}{\ssnl}
\ssnl 
\keepcomment{References to the literature.}{}

\nl
\bf Dual pairs of non-degenerate algebras \rm
\nl
\keepcomment{See reserve in artikel2ra.tex}
In the previous two items, we have obtained various algebras  of elements in the dual space $A'$. We have the pairing of $A$ with $A'$ and we can now restrict this to these various algebras. 

\ssnl
We now  investigate the properties of these pairings. 
\snl
For this recall the following notion found at different places in the literature, see e.g.\ Definition 3.1 in \cite{T-VD2} and Section 1 in \cite{La-VD1}.
\defin\label{defin:2.27e}
Let $A$ and $B$ be non-degenerate algebras and $(a,b)\mapsto \langle a,b\rangle$ a non-degenerate bilinear form on $A\times B$. 
Assume that there exist unital left and right actions of one algebra on the other defined by the equalities
\begin{align*}
\langle a,bb' \rangle &= \langle a\tl b,b' \rangle  &\langle a,b'b\rangle &=\langle  b\tr a, b'\rangle\\
\langle a'a,b \rangle &= \langle a',a\tr b\rangle &\langle aa',b \rangle &=\langle  a', b\tl a\rangle
\end{align*}
for $a,a'\in A$ and $b,b'\in B$. We say that the pairing admits actions.
\edefin

We have some general properties of these actions, see e.g.\ Section 1 of \cite{La-VD1}.

\prop
i) All these actions are faithful.
\ssnl
ii) If the algebra $A$ is idempotent, then the actions of $A$ on $B$ are non-degenerate while if $B$ is idempotent, the actions of $B$ on $A$ are non-degenerate.
\eprop
\bew
i) Let $a\in A$ and assume that $a\tr b=0$ for all $b\in B$. Then $\langle a'a,b\rangle =0$ for all $b$ and so, because the pairing is non-degenerate, $a'a=0$. This holds for all $a'$ and because the product is non-degenerate, we have $a=0$. This proves that the left action of $A$ on $B$ is faithful. Similarly for the three other actions.
\ssnl
ii) Let $b\in B$ and assume that $a\tr b=0$ for all $a\in A$. Then $\langle a'a,b\rangle=0$ for all $a,a'$. If $A$ is idempotent this implies that $b=0$. This means that the left action of $A$ on $B$ is non-degenerate. Similarly for the right action of $A$ on $B$ and for the actions of $B$ on $A$ when $B$ is idempotent. 
\ebew

There is still another property of actions. The left action of $A$ on $B$ is called unital if $B$ is spanned by elements of the form $a\tr b$ with $a\in A$ and $b\in B$. In general, we can not say anything about this.
\ssnl
We now look at our pairings of $A$ with the various subalgebras of $A'$ we have considered. These are first the smaller algebras $B^0_\ell$, $B^0_r$ and $B_0$. Secondly, we have the bigger algebras $B_\ell$, $B_r$ and $B$. 
\ssnl
All these algebras are non-degenerate. We refer to the Propositions  \ref{stel:2.6}, \ref{stel:2.7} and \ref{prop:2.11d} for the smaller ones and to Proposition \ref{prop:2.24d} and \ref{stel:2.25c} for the bigger ones.
\ssnl 
Moreover, all these pairings are still non-degenerate. Indeed, if $a\in A$ and $f(cac')=0$ for all $f\in A'$ and $c,c'\in A$ we must have that $a=0$ because the product in $A$ is non-degenerate. As all the algebras we consider contain all elements of the form $f(c\,\cdot\,c')$ and  this proves one side. For the other side, just observe that given $f\in A'$, then $f(a)=0$ for all $a$  implies $f=0$ by definition.
\nl
\bf The existence of the actions \rm
\nl
We first consider the actions of $A$ on these algebras.  

\prop\label{prop:2.29f} 
For all $a\in A$ and $b\in A'$ we have $b\tl a\in B^0_\ell$ and $a\tr b\in B^0_r$. These two actions are unital.
\eprop
\bew
i) Let $f$ be any linear functional on $A$ and $a\in A$. For all $x\in A$ we have $f(ax)=\langle x, b'\rangle$ where $b'=f(a\,\cdot\,)$. In other words we have $f\tl a=b'$. Because $B^0_\ell$ is the span of such alements, we have that the action is unital.
\ssnl
 ii) Let $f\in A'$ and $a\in A$. Then $f(xa)=\langle x,b'\rangle$ where $b'=f(\,\cdot\,a)$ Therefore $a\tr f=b'$. Again this action is unital.
\ebew

In particular, we can take for $b$ any element in one of the subspaces we consider and we always get $b\tl a\in B^0_\ell$ and $a\tr b\in B^0_r$ for all $a$. Further we  have that,  if $b\in B^0_\ell$, then $a\tr b\in B^0_\ell$  and hence also $a\tr b\in B_0$. Similarly, if $b\in B^0_r$ then $b\tl a\in B^0_r$ and hence also in $B_0$.
\ssnl
For the above results, we do not need the regularity of the canonical maps. When some of them are regular, we can say the following.

\prop\label {prop:2.30e}
If all canonical maps are regular, then $A'\tl A \subseteq B$ and $A\tr A'\subseteq B$.
\eprop

\bew
We always have $A' \tl A\subseteq B^0_\ell$. When all canonical maps are regular, we have $B_\ell\subseteq B$ by Proposition \ref{prop:2.15e}. Hence $A'\tl A \subseteq B$. Also, because $A\tr A'\subseteq B^0_r$ we have $A\tr A'\subseteq B$ by the same proposition.
\ebew

We see that in this case, $B$ is a two-sided $A$-module. It is not clear if this still holds for the smaller algebra $B_0$. One can also not expect that the actions of $A$ on $B$ are unital.
\nl
Next we consider the actions of the dual algebras on $A$.

\prop\label{prop:2.31f}
The left action of $B_r$ and the right action of $B_\ell$ on $A$ exist. In particular, the left and right actions of $B$ on $A$ exist.
\eprop

\bew
i) If $\omega_1\in A'$ and $\omega\in B_r$ we have, for all $a$, 
\begin{equation*}
(\omega_1\omega)(a)=\omega_1((\iota\ot\omega)\Delta(a)).
\end{equation*}
This means that the left action of $\omega$ on $a$ is given by $\omega\tr a=(\iota\ot\omega)\Delta(a)$.
\ssnl
ii) Similarly, if $\omega\in B_\ell$ and $\omega_2\in A'$ we have, for all $a$,
\begin{equation*}
(\omega\omega_2)(a)=\omega_2((\omega\ot\iota)\Delta(a)).
\end{equation*}
This means that the right action of $\omega$ on $a$ is given by $a\tl \omega=(\omega\ot\iota)\Delta(a)$.
\ebew

In particular, the left and right actions of $B$ on $A$ exist. We have seen in Proposition \ref{prop:2.30e} that the left and right actions of $A$ on $B$ exist when all canonical maps are regular. Then we can summarize and get the following property of the pairing of $A$ with $B$.

\prop
If all canonical maps are regular, then the pairing of $A$ with $B$ is an admissible pairing in the sense of Definition \ref{defin:2.27e}
\eprop

In Proposition \ref{prop:2.15e} we have seen that $B^0_\ell$ and $B^0_r$ belong to $B$ when all canonical maps are regular. Then these two algebras have left and right actions on $A$. This can be seen also directly.
\ssnl
Recall also that, for the dual algebras to be non-degenerate, we need that $\Delta$ is full. This implies the following for the actions.

\prop
The right actions of $B_\ell$ and $B^0_\ell$ and the left actions of $B_r$ and $B^0_r$  on $A$ are unital.
\eprop
\bew
If $b=f(c\,\cdot\,)$ where $f\in A'$ and $c\in A$, we have $a\tl b=(f\ot\iota)((c\ot 1)\Delta(a))$ and when $\Delta$ is full, such elements span all of $A$. It follows that the left action of $B^0_\ell$ on $A$ is unital. Because $B^0_\ell\subseteq B_\ell$ also the left action of $B_\ell$ on $A$ is unital. Similarly, if $b=f(\,\cdot\,c)$ we have $b\tr a=(\iota\ot f)(\Delta(a)(1\ot c))$ and again, because $\Delta$ is full, such elements span all of $A$.
\ebew

In the next section, we will illustrate some of these results. In particular, see the examples with infinite matrix algebras. We also refer to Remark \ref{opm:3.19e}.


%
%

\section{\hspace{-17pt}. More examples and special cases} \label{s:examp}  

In this section, we will consider some examples and special cases to illustrate various notions and results from the previous sections.
\nl
\bf Two trivial cases \rm
\nl
Let $A$ be a non-degenerate algebra $A$ and  $\Delta:A\to M(A\ot A)$ a linear map. Then we can consider the canonical maps as in Definition \ref{defin:1.2} of Section \ref{s:copr}).
\ssnl
It is rather easy to give (trivial) examples of a linear map  $\Delta: A\to M(A\ot A)$ where some of the canonical maps are regular and others are not. 

\voorb\label{voorb:4.1}
i) Define $\Delta(a)=a\ot 1$ for $a\in A$ and where $1$ is the identity in $M(A)$.
For the canonical maps we find
\begin{align*}
T_1(a\ot b)&=a\ot b
\tussenen 
T_2(c\ot a)=ca\ot 1\\
T_3(a\ot b)&=a\ot b
\tussenen
T_4(c\ot a)=ac\ot 1
\end{align*}
for all $a,b,c\in A$. 
\ssnl
ii) We see that $T_1$ and $T_3$ are the identiy maps from $A\ot A$ to itself. In particular they are regular. On the other hand, if $A$ does not have an identiy, $T_2$ and $T_4$ are not regular.
\evoorb

The common notion of coassociativity as in Definition \ref{defin:coass} can not be formulated. The same is true for the notions as in Definition \ref{defin:1.11a}. For this example however, $\Delta$ is a homomorphism and because $T_1$ is regular, Definition \ref{defin:1.6} does apply. 

\prop
The map $\Delta$ given by $\Delta(a)=a\ot 1$ is coassociative in the sense of Definition \ref{defin:1.6}.
\eprop
\bew
For $a,p,q\in A$ we find for the left hand side of Equation (\ref{eqn:1.7a})  
\begin{equation*}
((\iota\ot\Delta)(\Delta(a)(1\ot p)))(1\ot 1\ot q)=
((\iota\ot\Delta)(a\ot p))(1\ot 1\ot q)=
a\ot p\ot q
\end{equation*}
while for the right hand side
\begin{align*}
((\Delta\ot\iota)(\Delta(a)(1\ot q)))(1\ot p\ot 1)
&=((\Delta\ot\iota)(a\ot q))(1\ot p\ot 1)\\
&=(a\ot 1\ot q)(1\ot p\ot q)=a\ot p\ot q
\end{align*}
because $\Delta(p)(1\ot q)=p\ot q$.  We indeed, as expected, get the same expressions. Therefore $\Delta$ is coassociative in the sense of Definition \ref{defin:1.6}.
\ebew

We also have the following.

\prop\label{prop:3.3d}
As before, consider $\Delta(a)=a\ot 1$ for all $a$. 
If $A$ is idempotent, then $\Delta$ is a non-degenerate homomorphism as in Definition \ref{defin:nd}. The extension obtained in Proposition \ref{prop:ext} is given by $\Delta(m)=m\ot 1$ for $m\in M(A)$. The ones obtained in Proposition \ref{prop:ext2} are given by 
\begin{equation*}
(\iota\ot \Delta)(m)=m\ot 1
\tussen
(\Delta\ot\iota)(m)=\zeta_{12}(1\ot m) 
\end{equation*}
where $\zeta_{12}$ is the flip map $p\ot q\ot r\mapsto q\ot p\ot r$, extended to the muliplier algebras. Now $\Delta$ is coassociative in the sense of Definition \ref{defin:coass2}.
\eprop
\bew
The first statements are easily to verify. To show coassociativity we have
\begin{align*}
(\Delta\ot\iota)\Delta(a)&=\Delta(a)\ot 1=a\ot 1\ot 1\\
(\iota\ot\Delta)\Delta(a)&=a\ot\Delta(1)=a\ot 1\ot 1.
\end{align*}

\vskip -15pt
\ebew

It is clear that the result we are using for this example is not optimal. Indeed, one does not need that $A$ is idempotent to extend these maps. On simply can define $\Delta(m)=m\ot 1$ for $m\in M(A)$ and do something similar for the maps $\Delta\ot\iota$ and $\iota\ot\Delta$ in order to obtain coassociativity
\begin{equation*}
(\Delta\ot\iota)\Delta(a)=(\iota\ot\Delta)\Delta(a)=a\ot 1\ot 1.
\end{equation*}
We have included the above just to illustrate the general results in this case. See also Remark \ref{opm:1.20}. The point is that when $A$ is not idempotent, these extensions may not be the only ones.
\ssnl
Before we pass to some related examples, observe that we only have a single sided notion of a counit and fullness of the coproduct. 

\prop For this example, a left counit only exists if $A$ is trivial.
\eprop
\bew
Assume that $\varepsilon$ is a linear functional such that $(\varepsilon\ot\iota)(\Delta(a)(1\ot c))=ac$ for all $a,c$. Then $\varepsilon(a)c=ac$. In particular, $\varepsilon$ can not be $0$ because this would imply that the product in $A$ is triviallly $0$. Now take any $a$ such that $\varepsilon(a)=1$. Then $ac=c$ for all $c$. This implies that $a$ is a left unit. If we multiply with any $b$ and cancel $c$ we find that also $ba=b$ for all $b$. Hence $a$ is an identity. So $A$ is unital. We also get that $a=\varepsilon(a)1$ for all $a$. Hence $A$ is trivial. 
\ebew

Also observe that the left leg of $\Delta$ is all of $A$ but that we cannot define the right leg of $\Delta$ in $A$.
\ssnl
Here is a derived example.

\voorb\label{voorb:4.2}
Let $B$ and $C$ be two algebras and take $A=B\ot C$. Now define $\Delta$ on $A$ by $\Delta(b\ot c)=b\ot 1_C\ot 1_B\ot c$ where we have used $1_B$ for the identity in $M(B)$ and $1_C$ for the one in $M(C)$. We again have a homomorphism. In this case, if neither $B$ nor $C$ has an identity, none of the canonical maps will be regular. If the algebras are idempotent,  we have a non-degenerate coproduct and $\Delta$ is coassociative in the sense of Definition \ref{defin:coass2}.
\evoorb

We leave the details to the reader. Again it is not really necessary to have idempotent algebras.
\nl

We will now illustrate some of the notions and results from Section \ref{s:dual} for these examples.

\voorb
Consider the trivial example where $\Delta(a)=a\ot 1$ for all $a$. Because we have regularity of the map $T_1$ we can define $(\omega\ot\iota)\Delta(a)$ as a left multiplier. We get $\omega(a)1$. This will not be in $A$ if $A$ does not have an identity. Still, it is instructive to consider the case here with $1\in A$. For the product, dual to the coproduct, we find $(\omega_1\omega_2)(a)=\omega_1(a)\omega_2(1)$. This makes all of $A'$ into an associative algebra. However, the product is degenerate. If $\omega_2(1)=0$, then $\omega_1\omega_2=0$ for all $\omega_1$. This does not imply that $\omega_2=0$ (except when $A=\mathbb C 1$).
\evoorb

The case of Example \ref{voorb:4.2} is very similar and still not very interesting. 
\nl
\bf The infinite matrix algebra \rm
\nl
We now treat some examples involving infinite matrices.

\notat
Let $A$ be the algebra spanned by an infinite set of matrix elements $\{e_{ij}\mid i,j=1, 2, 3 \dots \}$.
\enotat
 Recall that $e_{ij}e_{kl}=\delta(j,k)e_{il}$ where $\delta$ is the Kronecker symbol.
 
 \prop\label{prop:3.7a}
 Fix $p,q$ and define $E_n$ in $A\ot A$ for any $n$ by 
 \begin{equation*}
E_n=\sum_{j=1}^n  e_{pj}\ot e_{jq}.
\end{equation*}
Then there is a multiplier $E\in M(A\ot A)$ such that for each $x\in A\ot A$ we have $E_nx=Ex$ and $xE_n=xE$ for $n $ large enough.
\eprop

 \bew
Given $p,q$,  we would like to define a multiplier $E$ by 
 \begin{equation*}
 E(e_{rs}\ot 1)=e_{ps}\ot e_{rq} \tussenen
 (1\ot e_{rs})E=e_{ps}\ot e_{rq}.
\end{equation*}
for all $r,s$. We have , for all $r,s$ and $r',s'$, 
\begin{align*}
(1\ot e_{r's'})(E(e_{rs}\ot 1))&=e_{ps}\ot e_{r's'}e_{rq}=\delta(r,s')e_{ps}\ot e_{r'q} \\
((1\ot e_{r's'})E)(e_{rs}\ot 1)&=e_{ps'}e_{rs}\ot e_{r'q}=\delta(r,s') e_{ps}\ot e_{r'q}.
\end{align*}
This is sufficient to prove that we indeed can define the multiplier $E$ of $A\ot A$. It is easy to see that 
$$E_n x\to Ex \tussenen xE_n\to xE$$ 
for all $x\in A\ot A$. In other words, given $x$ we have an $n_0$ such that $E_nx=x=xE_n$ when $n\geq n_0$.
 \ebew
 
The result seems to be obvious, but still, to be precise, we have to do some verifications as above.
\ssnl
We will use $\sum_j e_{pj}\ot e_{jq}$ to denote this multiplier.

\defin\label{defin:3.8d}
 We define the linear map $\Delta: A\to M(A\ot A)$ by 
 \begin{equation*}
\Delta(e_{pq})=\sum_j e_{pj}\ot e_{jq}.
\end{equation*}
\edefin

First we look for the canonical maps for this example.
\prop\label{prop:3.9d}
For this example, the canonical maps $T_1$ and $T_2$ are not  regular but the maps $T_3$ and $T_4$ are  regular.
 \eprop
 \bew
 i) The regularity of the maps $T_3$ and $T_4$ is a consequence of the formulas we had in the previous proof.
 \ssnl
 ii) On the other hand we have, for all $p,q$ and $r,s$, 
 \begin{align*}
T_1(e_{pq}\ot e_{rs})&=\sum_j e_{pj} \ot e_{jq}e_{rs}= \delta(q,r)\sum_j e_{pj} \ot e_{js}  \\
T_2(e_{pq}\ot e_{rs})&=\sum_j e_{pq}e_{rj}\ot e_{js}=\delta(q,r)\sum_je_{pj}\ot e_{js}.
\end{align*}
When $q=r$ these images do not belong to $A\ot A$.
\ebew
 
 We can use the  form of coassociativity with the maps $T_3$ and $T_4$ as in item ii) of Definition \ref{defin:coass}. However, we want a stronger version of this property. We will prove it later (see Proposition \ref{prop:1.6} below), first we show the following.
 
\prop
There is a counit, defined by $\varepsilon(e_{pq})=\delta(p,q)$. 
\eprop
\bew
i) We have
\begin{align*}
(\varepsilon\ot\iota)((1\ot e_{kl})\Delta(e_{pq}))
&=(\varepsilon\ot\iota)(e_{pl}\ot e_{kq})\\
&=\delta(p,l) e_{kq}=e_{kl}e_{pq}
\end{align*}
Similarly we get $(\iota\ot\varepsilon)(\Delta(e_{pq})(e_{kl}\ot 1))=e_{pq}e_{kl}$.
\ebew 

Remark that we use the maps $T_3$ and $T_4$ to characterize the counit for this example. See a remark following Definition \ref{defin:1.15c}.
\ssnl
It is not a homomorphism. Indeed, we have
\begin{equation*}
\varepsilon(e_{pq}e_{kl})=\delta(q,k)\varepsilon(e_{pl})=\delta(q,k)\delta(p,l)
\end{equation*}
and e.g. $\varepsilon(e_{12}e_{21})=1\neq 0= \varepsilon(e_{12})\varepsilon(e_{21})$.  
\nl 
\bf The algebras $B_\ell$ and $B_r$ for this example. \rm
\nl
 We have the following characterization of the spaces $B_\ell$ and $B_r$.
 
\prop
i) An element $\omega$ of $A'$ belongs to $B_\ell$ if and only if $j\mapsto \omega(e_{pj})$ has finite support for each $p$. It belongs to $B_r$ if and only if $j\mapsto \omega(e_{jq})$ has finite support for each $q$.
\ssnl
ii) If either $\omega_1\in B_\ell$ and $\omega_2\in A'$ or $\omega_1 \in A'$ and $\omega_2\in B_r$, the product $\omega_1\omega_2$ satisfies 
\begin{equation*}
(\omega_1\omega_2)(e_{pq})=\sum_j \omega_1(e_{pj})\omega_2(e_{jq}).
\end{equation*}
\eprop
\bew
i) Let $\omega\in A'$. Then we have
\begin{equation*}
(\omega\ot\iota)\Delta(e_{pq})
=\sum_j \omega(e_{pj})\,e_{jq}.
\end{equation*}
For $(\omega\ot \iota)\Delta(e_{pq})$ to be in $A$ we need a finite sum here. In other words, $\omega\in B_\ell$ if and only if $j\mapsto \omega(e_{pj})$ has finite finite support for each $p$. Similarly $\omega\in B_r$ if and only if $j\mapsto \omega(e_{jq})$ has finite support for each $q$.
\ssnl
ii) Now take $\omega_1\in B_\ell$ and $\omega_2\in A'$. Then we get, for the product as defined in Definition \ref{defin:2.7a},
\begin{equation*}
(\omega_1\omega_2)(e_{pq})=\sum_j \omega_1(e_{pj})\omega_2 (e_{jq}).
\end{equation*}
The sum is well-defined because $\omega_1\in B_\ell$ so that there are only finitely many indices $j$ for which $\omega_1(e_{pj})\neq 0$. Similarly weh $\omega_1\in A'$ and $\omega_2\in B_r$. Here we use item ii) of Definition \ref{defin:2.7a} and that $j\mapsto \omega_2(e_{jq})$ has finite support.
\ebew

We can verify that the product of elements in $B_\ell$ is again in $B_\ell$. Indeed, because $\omega_1$ is in $B_\ell$ we have, given $p$ only finitely many $j$ so that $\omega_1(e_{pj})$ is non-zero. For this set of indices $j$ we have only finitely many $q$ such that $\omega_2(e_{jq})$ is non-zero. Hence, given $p$ there are only finitely many $q$ such that $(\omega_1\omega_2)(e_{pq})$ is non-zero. Therefore $\omega_1\omega_2\in B_\ell$. 
\ssnl
In Definition \ref{defin:2.8} we formulated the stronger form of coassociativity in the case where $T_1$ and $T_2$ are regular. Here however we have that  $T_3$ and $T_4$ are regular. We therefore modify this form of coassociativity.  Because $T_3$ and $T_4$ are regular we consider item ii) of Definition \ref{defin:coass}.
 Given $a,b,c\in A$ we have
 \begin{equation*}
(1\ot 1\ot c)((\iota\ot\Delta)(\Delta(a)(b\ot 1)=((\Delta\ot\iota)((1\ot c)\Delta(a)))(b\ot1\ot 1).
\end{equation*}
Now we take two linear functionals $f$ and $g$ and we let $\omega_1=f(\,\cdot\,b)$ and $\omega_2=g(c\,\cdot\,)$. When we apply $f\ot \iota\ot g$ to the above equation we get
\begin{equation*}
(\iota\ot \omega_2)\Delta((\omega_1\ot\iota)\Delta(a))
=(\omega_1\ot\iota)\Delta((\iota\ot\omega_2)\Delta(a)).
\end{equation*}
In this way, we get a stronger form of coassociativity as in Definition \ref{defin:2.8} but now for the case where $T_3$ and $T_4$ are regular. 
\keepcomment{\ssnl We may have to formulate this a bit better?}{}
\ssnl
We now show that this holds for our example.

\prop \label{prop:1.6}
For all $\omega_1\in B_\ell$ and $\omega_2\in B_r$ we have
\begin{equation*}
(\iota\ot \omega_2)\Delta((\omega_1\ot\iota)\Delta(a))
=(\omega_1\ot\iota)\Delta((\iota\ot\omega_2)\Delta(a))
\end{equation*}
for all $a$.
\eprop
\bew
For all $p,q$ we have
\begin{align*}
(\iota\ot \omega_2)\Delta((\omega_1\ot\iota)\Delta(e_{pq}))
&=\sum_j (\iota\ot \omega_2)\Delta(\omega_1(e_{pj})e_{jq}) \\
&=\sum_{j,k} \omega_1(e_{pj})\, e_{jk} \,\omega_2(e_{kq}).
\end{align*}
The sum over $j$ is a finite sum because $\omega_1\in B_\ell$ while the one over $k$ is a finite sum because $\omega_2\in B_r$. Similarly we get
\begin{align*}
(\omega_1\ot\iota)\Delta((\iota\ot \omega_2)\Delta(e_{pq}))
&=\sum_k (\omega_1\ot \iota)\Delta((e_{pk}\omega_2(e_{kq})) \\
&=\sum_{j,k} \omega_1(e_{pj})\, e_{jk}\, \omega_2(e_{kq}).
\end{align*}

\vskip-10pt\ebew

Formally, this is easier:
\begin{align*}
(\iota\ot\Delta)\Delta(e_{pq})&=\sum_k e_{pk} \ot\Delta(e_{kq}) = \sum_{j,k} e_{pj}\ot e_{jk}\ot e_{kq}\\
(\Delta\ot\iota)\Delta(e_{pq})&=\sum_j \Delta(e_{pj}) \ot e_{jq} = \sum_{j,k} e_{pj}\ot e_{jk}\ot e_{kq}
\end{align*}

We now obtain more results about the dual algebras 
$B_\ell$ and $B_r$ and their intersection $B$. We need the following result.

\prop\label{prop:3.14e}
We have
\begin{equation*}
\omega_2((\omega_1\ot\iota)\Delta(a))
=\omega_1((\iota\ot\omega_2)\Delta(a))
\end{equation*}
when $\omega_1\in B_\ell$ and $\omega_2\in B_r$.
\eprop

\bew
Let $a=e_{pq}$. For the left hand side we find
\begin{equation*}
\omega_2((\omega_1\ot\iota)\Delta(e_{pq}))=\omega_2(\textstyle\sum_j \omega_1(e_{pj})e_{jq})=\textstyle\sum_j \omega_1(e_{pj})\omega_2(e_{jq}).
\end{equation*}
The last equality is justified because $\omega_1\in B_\ell$ so that we have a finite sum. Similarly, for the right hand side we get
\begin{equation*}
\omega_1((\iota\ot\omega_2)\Delta(e_{pq}))=\omega_1(\textstyle\sum_j \omega_2(e_{jq})e_{pj})=\textstyle\sum_j \omega_1(e_{pj})\omega_2(e_{jq}).
\end{equation*}
Now, the last equality is true because $\omega_2\in B_r$ so that we have a finite sum.
\ebew

For the proof of this result in the general case, as in Proposition \ref{prop:2.10b} we used the counit and that it is a homomorphism. for this example however, the counit is not a homomorphism and it is a bit strange that still the result here is true.
\ssnl
As a consequence, we have the result of Proposition \ref{stel:2.25c} for this case and the intersection $B$ of $B_\ell$ and $B_r$ is an algebra for the product inherited from the two.

\prop
The counit belongs to $B$ and it is a unit for alle these algebras.
\eprop
\bew
We have seen in the proof of Proposition \ref{prop:3.9d} that
\begin{equation*}
(\varepsilon\ot\iota)\Delta(e_{pq})=e_{pq}
\tussenen
(\iota\ot\varepsilon)\Delta(e_{pq})=e_{pq}
\end{equation*}
for all $p ,q$. This proves that $\varepsilon\in B_\ell$ and $\varepsilon\in B_r$. To prove that it is a unit, just apply any $\omega$ of $A'$ to these equations to obtain that
\begin{equation*}
\varepsilon\omega(e_{pq})= \omega(e_{pq})
\tussenen
\omega\varepsilon(e_{pq})= \omega(e_{pq})  
\end{equation*}
for all $p,q$.
\ebew

In fact we also have the following.

\prop
For all $\omega_1\in B_\ell$ and $\omega_2\in B_r$ we have, for all $a$,
\begin{equation*}
\varepsilon((\omega_1\ot\iota)\Delta(a))=\omega_1(a)
\tussenen
\varepsilon((\iota\ot\omega_2)\Delta(a))=\omega_2(a).
\end{equation*}
\eprop
\bew
\begin{equation*}
\varepsilon((\omega_1\ot\iota)\Delta(e_{pq}))=\sum_j \omega_1(e_{pj})\varepsilon(e_{jq})=\omega_1(e_{pq}).
\end{equation*}
Similarly $\varepsilon((\iota\ot\omega_2)\Delta(e_{pq}))=\omega_2(e_{pq})$.
\ebew

As mentioned earliier (see Remark \ref{opm:2.17c}) the above result also follows from Proposition \ref{prop:3.14e}.

\nl 
\bf The algebras $B^0_\ell$, $B^0_r$ and $B_0$ \rm
\nl
In this case, the algebra $B_\ell$ is strictly bigger than the space $B^0_\ell$ spanned by elements in $A'$ of the form $f(c\,\cdot\,)$ defined in \ref{defin:2.1d}. Similarly, the algebra $B_r$ is stricly bigger than the space $B^0_r$ spanned by elements in $A'$ of the form $f(\,\cdot\,c)$, defined in Definition \ref{defin:2.1d}.
We show this in the next proposition. 

\prop
The space $B^0_\ell$ consists of linear functionals $\omega$ on $A$ with the property that the map $j\mapsto \omega(e_{pj})$ has finite support for all $p$, independent of $p$. It is a subalgebra of $B_\ell$. Similarly, the space $B^0_r$  
are the linear functionals $\omega$ such $j\mapsto \omega(e_{jq})$ has finite support for all $q$,  indepentent of $q$. It is a subalgebra of $B_r$.
\eprop
\bew
i) Let $\omega$ be a finite linear combination of linear functionals of the form $f(e_{rs}\,\cdot\,)$. Denote by $I_\omega$ the set of indices  $s$ that appear in this form. Then $\omega(e_{pq})=0$ for all $q$ when $p\notin I_\omega$. For two such linear functionals $\omega_1, \omega_2$ we have $(\omega_1\omega_2)(e_{pq})=\sum_j \omega_1(e_{pj})\omega_2(e_{jq})$. For $p\notin I_{\omega_1}$ we have $\omega_1(e_{pj})=0$ for all $j$. Then also $(\omega_1\omega_2)(e_{pq})=0$ for all $q$. This proves that also $\omega_1\omega_2$ belongs to this space.
\ssnl
ii) Similarly if $\omega$ is a linear combination of linear functionals of the form $f(\,\cdot\,e_{rs})$. Now we take for $J_\omega$ the set of indices $r$ that appear in this expression.
\ebew

It is clear that the intersection $B_0$  of these two algebras is the algebra of functionals $\omega$ with the property that $\omega(e_{pq})=0$ except for finitely many pairs of indices. This algebra is isomorphic with the original algebra as we see in the next proposition.


\prop
Define elements $f_{rs}$ for all $r,s$  in $A'$ by
\begin{equation*}
\langle e_{pq},f_{rs} \rangle= \delta(p,r)\delta(q,s).
\end{equation*}
Then they are matrix elements and span  the algebra $B_0$.
\eprop
\bew
We clearly  have that these elements belong to $B_0$ and that the linear span of them is all of $B_0$. 
For all $p,q$, $r,s$ and $r',s'$ we have, using the definition of the product in the dual algebras, 
\begin{align*}
\langle e_{pq}, f_{rs}f_{r's'}\rangle
&=\sum_j \langle e_{pj},f_{rs}\rangle \langle e_{jq},f_{r's'} \rangle \\
&=\sum_j \delta(p,r)\delta (j,s) \delta (j,r')\delta(q,s') \\
&=\delta(p,r)\delta(s,r')\delta(q,s') \\
&=\delta(s,r')\langle e_{pq},f_{rs'}\rangle.
\end{align*}
We see that $f_{rs}f_{r's'}=\delta(s,r')f_{rs'}$.
\ebew

The algebra $B_0$ is isomorphic with the original algebra and so we get a self-pairing of $A$. This is the right place to formulate the following important remark.

\opm\label{opm:3.19e}
We could have treated this example from an other point of view. 
\ssnl
Indeed, consider two copies of the infinite matrix algebra. To be consistent with the previous approach, we denote these two algebras with $A$ and $B_0$. We use $(e_{pq})$ for the set of matrix units that spans $A$ and $(f_{pq})$ for the ones that span $B_0$.  Now we define a pairing of $A$ with $B_0$ by
\begin{equation*}
\langle e_{pq},f_{rs} \rangle= \delta(p,r)\delta(q,s).
\end{equation*}
One can easily verify that this is an admissible pairing as in Definition \ref{defin:2.27e}. One then defines a coproduct $\Delta:A\mapsto B\ot B)'$ by $\langle \Delta(a),b\ot b'\rangle = \langle a,bb'\rangle$. The next step is to show that actually $\Delta(A)\subseteq M(A\ot A)$ , where we have the extension of the pairing from $(A\ot A)\times (B_0\ot B_0)$ to $M(A \ot A)\times (B_0\ot B_0)$. 
\ssnl
All objects and results we had for this example, can be obtained from this starting point. 
\ssnl
For such an approach, we refer to \cite{La-VD1}.
\eopm

None of these algebras contains the counit. This is compatible with the fact that it is not a homomorphism because of Proposition \ref{prop:2.11c}. 
Also remark that in Proposition \ref{prop:2.5c}, it was needed to have that $\Delta$ is a homomorphism to obtain that these spaces are algebras. That is not the case here for this example.

\nl
As the algebras are non-degenerate, it makes sense to find the mulltiplier algebras. 
\keepcomment{\rood \ssnl We have to refer to our paper with Joost on the infinite matrix algebras.}{} 

\prop
The multiplier algebra of $B_0$ is $B$.
\eprop
\bew
First remark that any multiplier of $B_0$ is given by an element $\omega\in A'$. To see this, let $m$ be a multiplier of $B_0$ and look at $f_{pp}mf_{qq}$. This is a multiple of $f_{pq}$. We can define $\omega(e_{pq})$ by 
$$f_{pp}mf_{qq}=\omega(e_{pq})f_{pq}.$$ 
We must have $f_{pp}m=\sum_q \omega(e_{pq}) f_{pq}$ and because this belongs to $B_0$ we must have, for each $p$ only finitely many $q$ with $\omega(e_{pq})\neq 0$. Similarly $m f_{qq}=\sum_p \omega(e_{pq}) f_{pq}$ and we must have, for each $q$ only finitely many $p$ with $\omega(e_{pq})\neq 0$. Hence $\omega\in B$.
\ebew
In a similar way, we get the following.
\prop
The multiplier algebra of $B^0_\ell$ is $B_\ell$. 
\eprop
\bew
i) Because the matrix elements $f_{pq}$ belong to $B^0_\ell$ we have again that any multiplier of $B^0_\ell$ is given by a linear functional $\omega$ on $A$. For any $\omega_1\in B^0_\ell$ we have 
\begin{equation*}
(\omega_1\omega)(e_{pq})=\sum_j \omega_1(e_{pj})\omega(e_{jq}).
\end{equation*}
Because $\omega_1\in B^0_\ell$ there is a finite subset of indices $J$ so that $\omega_1(e_{pj})=0$ for all $p$ and all $j\notin J$. Then we will have similar result for the product if we have for each $j\in J$ that $\omega(e_{jq})=0$ except for finitely many $q$. This is the case if and only if $\omega\in B_\ell$. 
\ssnl
ii) On the other hand, assume that $\omega_1\in B_\ell$ and $\omega_2\in B^0_\ell$. Then there is a finite set of indices $I$ such that $\omega_2(e_{jq})=0$ for all $j$ when $q\notin I$. Then $(\omega_1\omega_2) (e_{pq})=0$ for all $p$ when $q\notin I$. Hence again $\omega_1\omega_2\in B^0_\ell$.
\ssnl
iii) Together we find that the multiplier algebra of $B^0_\ell$ is precisely  $B_\ell$.
\ebew

Similarly, the multiplier algebra of $B^0_r$ is $B_r$.

\nl
\bf Properties of the pairing between $A$ and these algebras\rm
\nl
We have pairings of $A$ with any of the algebras $B_\ell$, $B_r$ and $B$ and we further use the pairing notation.
\ssnl
What about the actions? We can look at the general results, but we have to take into account that only the canonical maps $T_3$ and $T_4$ are regular for this case.
\ssnl
First we consider the left and the right actions of $A$. For these actions, the regularity of the canonical maps does not play a role. So we have the result as formulated in Proposition \ref{prop:2.29f}.

\prop
For all $a\in A$ and $b\in A'$ we have $b\tl a\in B^0_\ell$ and $a\tr b\in B^0_r$. These two actions are unital.
\eprop

We also have that $a\tr b\in B_0$ if $b\in B^0_\ell$ and $b\tl a\in B_0$ if $b\in B^0_\ell$.
\ssnl
Because the canonical maps $T_3$ and $T_4$ are regular, it follows from the results in Proposition \ref{prop:2.15e} that $A'\tl A\subseteq B_r$ and $A\tr A'\subseteq B_\ell$. In particular, $B_r$ is a right $A$-module and $B_\ell$ is a left $A$-module. Here we can not apply Proposition \ref{prop:2.30e} because the maps $T_1$ and $T_2$ are not regular.
\ssnl
Next we look at the actions of the dual algebra on $A$. We have the result of Proposition \ref{prop:2.31f}. 

\prop
The left action of $B_r$ and the right action of $B_\ell$ on $A$ exist. In particular, the left and right actions of $B$ on $A$ exist.
\eprop

The result remains true as we can still define the spaces $B_\ell$ and $B_r$ with the regularity of the maps $T_3$ and $T_4$, see a remark following Definition \ref{defin:2.6}.
\ssnl
For the smaller algebras we get the following.

\prop We have
\begin{equation*}
A\tl B^0_\ell=A
\tussenen
B^0_r\tr A= A.
\end{equation*}
\eprop

\bew
We have $e_{pq}\tl b=\sum_j \langle e_{pj},b \rangle e_{jq}$ and if $b$ is the functional in $A'$ that is $0$ on all elements except on $e_{pp}$, we get  $e_{pq}\tl b=e_{pq}$ when $\langle e_{pp},b\rangle=1$. This proves that $A\tl B^0_\ell=A$. Similarly $B^0_r\tr A= A$.
\ebew

In fact we even get 
\begin{equation*}
A\tl B_0=A
\tussenen
B_0\tr A= A.
\end{equation*}
where $B^0=B^0_\ell\cap B^0_r$.
\keepcomment{ We need to summarize and compare with the general results.}{}
\ssnl
Because we also have that $A\tr B_0\subseteq B_0$ and $B_0\tl A\subseteq B_0$, 
we conclude that we have an admissible pair of $A$ with $B_0$. 

\nl
\bf More examples with infinite matrix algebras \rm
\nl
In what follows we will use $C$ for the infinite matrix algebra, spanned by an infinite set of matrix elements $(e_{pq})$ where $p,q=1, 2, 3, \dots$ as in the previous item. 

\voorb\label{voorb:4.3}
i)
Let $A=C\ot C$ and define $\Delta$ on $A$ by $\Delta(c_1\ot c_2 )=c_1\ot E\ot c_2$ where $E=\sum_j e_{1j} \ot e_{j1}$. This element $E$ is well-defined in the multiplier algebra $M(C\ot C)$ as we have seen in Proposition \ref{prop:3.7a}.
It follows that $\Delta(a)$ belongs to the multiplier algebra $M(A\ot A)$.
\ssnl
iii) As we have seen in the proof of Proposition \ref{prop:3.7a} we have
\begin{equation}
E(e_{rs}\ot 1)=e_{1s}\ot e_{r1} 
\tussenen
 (1\ot e_{rs})E=e_{1s}\ot e_{r1}.   \label{eqn:3.1}                                    
\end{equation}
We can conclude from the equalities above also that $T_3$ and $T_4$ are regular. On the other hand
\begin{equation*}
E(1\ot e_{11})=E
\tussenen
(e_{11}\ot 1)E=E
\end{equation*}
and consequently, the canonical maps $T_1$ and $T_2$ will not be regular. 
\ssnl
iv) One can verify that 
$\Delta$ is coassociative in the sense of Definition \ref{defin:coass}. Essentially we get
\begin{align}
(\Delta\ot\iota)\Delta(c\ot c')&=c\ot E\ot E\ot c'  \label{eqn:3.2a}\\
(\iota\ot\Delta)\Delta(c\ot c')&=c\ot E\ot E\ot c' \label{eqn:3.2b}
\end{align}
for $c,c'\in C$. To make these formulas precise, we should multiply the first equation on the left with $e_{rs}$ in the fifth factor, and use the second formula in Equation (\ref{eqn:3.1}). Also we have to multiply the second equation on the right with $e_{r's'}$ in the second factor, and use the first formula in Equation (\ref{eqn:3.1}).
\evoorb

We now look for the dual algebras.

\prop
A linear functional $\omega$ on $A$ belongs to $B_\ell$ if and only if $n\mapsto \omega(c\ot e_{1n})$ has finite support for all $c\in C$. Similarly, it belongs to $B_r$ if and only if $n\mapsto \omega(e_{n1}\ot c)$ has finite support for all $c\in C$.
\eprop

\bew

For any linear functional $\omega$ on $A$ we get
\begin{equation*}
(\omega\ot\iota)\Delta(c\ot c')=\sum_n \omega(c \ot e_{1n})\, e_{n1}\ot c'.
\end{equation*}
This will be an element of $A$ if $\omega(c\ot e_{1n})$ is $0$ except for finitely many $n$. Similarly we have
\begin{equation*}
(\iota\ot\omega)\Delta(c\ot c')=\sum_n \omega(e_{n1}\ot c')\,c\ot e_{1n}.
\end{equation*}
This will be an element in $A$ if $\omega(e_{n1}\ot c')=0$ except for finitely many $n$.
\ebew

For the product $\omega_1\omega_2$ we get
\begin{equation*}
(\omega_1\omega_2)(c\ot c')=\sum_n \omega_1(c\ot e_{1n})\,\omega_2(e_{n1}\ot c').
\end{equation*}
We see that we get a finite sum if either $\omega_1\in B_\ell$ or $\omega_2\in B_r$.
\ssnl
The coproduct is not full and indeed, this product is degenerate. If $\omega_1(c\ot e_{1n})=0$ for all $n$, we get $\omega_1\omega_2=0$ for all $\omega_2$ while if $\omega_2(e_{n1}\ot c)=0$ for all $n$, we have $\omega_1\omega_2=0$ for all $\omega_1$.
\ssnl
It is also possible to characterize the smaller algebras $B^0_\ell$, $B^0_r$ and $B_0$. They are strictly smaller than the larger ones $B_\ell$, $B_r$ and $B$. 
\keepcomment{\rood Check and look for details. Do we include them or leave them as an exercise for the reader?}{}
\ssnl
The following example is similar to the previous one, but does not involve the infinite matrix algebra.

\voorb\label{voorb:3.32e}
Let $X$ be any set and $P$ the algebra $F(X)$ of complex functions with finite support on $X$, endowed with pointwise operations. Take for $A$ the algebra $P\ot P$. It is naturally identified with the algebra $F(X\times X)$ of complex functions with finite support in $X\times X$. Denote by $\delta_x$ the function on $X$ with the value $1$ in the point $x$ and $0$ in all other points. Define a coproduct $\Delta$ on $A$ by $\Delta(p\ot p')=p\ot E\ot p'$ where $E=\sum_{x\in X}\delta_x\ot\delta_x$. 
It is easy to verify that $E$ is a multiplier of $C\ot C$ and that $E^2=E$. It follows that $\Delta$ is a homomorphism of $A$ in $M(A\ot A)$. 
\ssnl
All the canonical maps are regular. And $\Delta$ is coassociative because
\begin{equation*}
(\iota\ot\Delta)\Delta(p\ot p')=p\ot E\ot E\ot p'
\quad\text{and}\quad
(\Delta\ot\iota)\Delta(p\ot p')=p\ot E\ot E\ot p'.
\end{equation*}
\ssnl
For any linear function $\omega$ on $A$ we get
\begin{align*}
(\omega\ot\iota)\Delta(\delta_y\ot\delta_z)
&=\sum_x \omega(\delta_y\ot\delta_x)\,\delta_x\ot\delta_z\\
(\iota\ot\omega)\Delta(\delta_y\ot\delta_z)
&=\sum_x \omega(\delta_x\ot\delta_z)\,\delta_y\ot\delta_x.
\end{align*}
In both cases, we get elements in the multiplier algebra $M(A)$. The first one belongs to $A$ if and only if, for all $y$, $\omega(\delta_y\ot\delta_x)=0$ except for finitely many $x$. This defines $B_r$. The second on belongs to $A$ if and only if, for all $z$, $\omega(\delta_x\ot\delta_z)=0$ except for finitely many $x$. This defines $B_\ell$. For the two elements to belong to $A$ we simply need that the function $x\mapsto \omega(\delta_y\ot\delta_x)$ has finite support for each $q$ and that $x\mapsto \omega(\delta_x\ot\delta_z)$ has finite support for each $z$. Then we get $B_r\cap B_\ell$.
\ssnl
The product is given by
\begin{equation*}
(\omega_1\omega_2)(\delta_y\ot\delta_z)=\sum_x \omega_1(\delta_y\ot \delta_x)\omega_2(\delta_x\ot\delta_z).
\end{equation*}
This product is defined when either $\omega_1\in B_\ell$ belongs or $\omega_2\in B_r$. 
\evoorb

We now consider again the infinite matrix algebra $C$ spanned by the matrix elements $\{e_{ij}\}$. In what follows, we will denote by $P$ the abelian subalgebra of $C$ spanned by the elements $\{e_{ii}\}$. We write $p_j$ for $e_{jj}$.

\voorb
i) Let $A=P\ot C$ and define $\Delta$ on $A$ by $\Delta(p\ot c)=p\ot E\ot c$ where $E=\sum_, e_{1j} \ot p_j$. This infinite sum is well-defined in the multiplier algebra $M(C\ot  P)$. Indeed,  we have
\begin{equation}
E(1\ot p_r)=e_{1r}\ot p_r
\tussenen
(1\ot p_r)E=e_{1r}\ot p_r. \label{eqn:3.4a}
\end{equation}
As a consequence we get that $\Delta(p\ot c)$ belongs to the multiplier algebra of $A\ot A$.
\ssnl
ii) From the Equations \ref{eqn:3.4a}, it follows that $T_1$ and $T_3$ are regular. 
Indeed
\begin{align*}
\Delta(p\ot c)(1\ot 1\ot p_r\ot c')&=(p\ot E\ot c)(1\ot 1\ot p_r\ot c')=p\ot e_{1r}\ot p_r\ot cc'\\
(1\ot 1\ot p_r\ot c')\Delta(p\ot c)&=(1\ot 1\ot p_r\ot c')(p\ot E\ot c)=p\ot e_{1r}\ot p_r\ot c'c
\end{align*}
for all $p\in P$ and $c,c'\in C$. Also $T_4$ is regular, but $T_2$ is not. This follows from 
\begin{equation*}
E(e_{rs}\ot 1)=e_{1s}\ot e_{rr} 
\tussenen
(e_{11}\ot 1)E=E.
\end{equation*}
We can formulate coassociativity as in Definition \ref{defin:coass}. We get (formally)
\begin{equation*}
(\Delta\ot\iota)\Delta(p\ot c)=p\ot E\ot E\ot c
\end{equation*}
for all $p\in P$ and $c\in C$. 
\evoorb

Still, for this example, $\Delta$ is not a homomorphism. It is also not full. We will now modify this example, so as to get full  coproducts of this type that are homomorphisms.
\ssnl
First, in the next example, we modify the previous one so that $\Delta$ is a homomorphism.

\voorb\label{voorb:4.4}
i) 
Assume that $q_j$ is an idempotent in $C$ for all $j$. We can define $E=\sum_j q_j\ot p_j$ in the multiplier algebra of $C\ot P$ and we have 
\begin{equation*}
E(1\ot p_r)=(1\ot p_r)E=q_r\ot p_r
\end{equation*}
for all $r$. Because $q_r$ is an idempotent for all $r$, we have that $E^2=E$.
\ssnl
ii) As in the previous example, we let $A$ be the algebra $P\ot C$ and we define $\Delta$ on $A$ by
$\Delta(p\ot c)=p\ot E\ot c$.
As in the previous example, this is a well-defined linear map from $A$ to the multiplier algebra $M(A\ot A)$. In this case, because $E$ is an idempotent, $\Delta$ is a homomorphism.
\ssnl
iii) For this coproduct, the canonical maps $T_1$ and $T_3$ are regular. Indeed
\begin{align*}
\Delta(p\ot c)(1\ot 1\ot p_r\ot c')&=(p\ot E\ot c)(1\ot 1\ot p_r\ot c')=p\ot q_r\ot p_r\ot cc'\\
(1\ot 1\ot p_r\ot c')\Delta(p\ot c)&=(1\ot 1\ot p_r\ot c')(p\ot E\ot c)=p\ot q_r\ot p_r\ot cc'
\end{align*}
for all $p\in P$ and $c,c'\in C$. Because $\Delta$ is a homomorphism, we can formulate coassociativity as in Definition \ref{defin:1.6}. We get (formally)
\begin{equation*}
(\Delta\ot\iota)\Delta(p\ot c)=p\ot E\ot E\ot c
\end{equation*}
for all $p\in P$ and $c\in C$. 
\evoorb

In general we can not say anything about the regularity of the maps $T_2$ and $T_4$. Also, for the above example,  the coproduct is not necessarily full. The right leg of $E$ is all of $P$ so that the right leg of $\Delta$ is all of $A$. However, the left leg of $E$ is spanned by the elements $q_r$ and we have no information about this.
\ssnl
We pass to more specific cases by making appropriate choices for the idempotents $q_r$.
\ssnl
We could take $q_r=p_r$ for all $r$. Then the coproduct is a regular homomorphism, but the left leg is only $P$. We will consider this example later (see Example \ref{voorb:3.32e}).
\ssnl
We take a more sophisticated choice in the following proposition. We take $P$, $C$  and $A$ as in the previous proposition.

\prop
Let $q_1=e_{11}$ and $q_n=e_{n1}+e_{nn}$ for $n>1$.  This is an idempotent for all $n$. Define again $E=\sum_j q_j\ot p_j$ and $\Delta$ as before with this element $E$. Then $T_1$, $T_2$ and $T_3$ are regular but $T_4$ is not.
\eprop
\bew
i) First remark that $q_n$ is also an idempotent for $n>1$. Indeed
\begin{align*}
q_nq_n&=(e_{n1}+e_{nn})(e_{n1}+e_{nn})\\
&=e_{n1}e_{n1}+e_{n1}e_{nn}+e_{nn}e_{n1}+e_{nn}e_{nn}\\
&=0 + 0 + e_{n1}+e_{nn}=q_n.
\end{align*}
\ssnl
ii) We know from the previous result that the maps $T_1$ and $T_3$ are regular. Further, for all $r,s$ we have,
\begin{equation*}
(e_{rs}\ot 1)E=e_{rs}e_{11}\ot p_1+\sum_{j>1}\,e_{rs}(e_{j1}+e_{jj}) \ot  p_j
\end{equation*}
For $s=1$ this gives $e_{r1}\ot p_1$ and for $s>1$ we get $(e_{r1}+e_{ss})\ot p_s$. It follows that $T_2$ is regular. On the other hand
\begin{equation*}
E(e_{11}\ot 1)= e_{11}\ot p_1 + \sum_{j>1}^\infty  e_{j1}\ot p_j
\end{equation*}
and we see that $T_4$ is not regular.
\ebew

For this example, the coproduct is not full. The left leg of $E$ is spanned by the elements $q_n$ and so the left leg of $\Delta$ is spanned by elements of the form $c\ot q_n$.
\ssnl
We now try to modify this example to get a full coproduct. 

\voorb\label{voorb:3.24}
Denote by $C_u$ the subalgebra of upper diagonal matrices in $C$. It is spanned by the matrix elements $\{e_{ij}\mid i\leq j \}$. Remark that this subalgebra has local units. Indeed $e_{ij}e_{jj}=e_{ii}e_{ij}=e_{ij}$. In particular, it is non-degenerate. It is non-unital. The algebra $P$ is still a subalgebra.
\ssnl
Let $E=\sum_{i\leq j} q_{ij}\ot (p_i\ot p_j)$ where
\begin{equation*}
q_{jj} = e_{jj} \tussenen q_{ij}=e_{ii}+te_{ij} \text{ if } i<j.
\end{equation*}
Remark that $q_{ij}$ is an idempontent for all pairs. Indeed, if $i<j$ we have
\begin{equation*}
q_{ij}^2=(e_{ii}+te_{ij})(e_{ii}+te_{ij})=e_{ii} +t e_{ij} + 0+0=q_{ij}.
\end{equation*}
The subalgebra $P_u$, spanned by the elements $\{p_i\ot p_j\mid i<j\}$ of $P\ot P$, is isomorphic with $P$. In this way, we get a multiplier in $C_u\ot P_u$. It is an idempotent.
\ssnl
What about the regularity properties. First we have $E(1\ot p_{ij})=(1\ot p_{ij})E=q_{ij}\ot p_{ij}$. We have used $p_{ij}$ for $p_i\ot p_j$. So $E(1\ot P_u)\subseteq C_u\ot P_u$ as well as $(1\ot P_u)E\subseteq C_u\ot P_u$. On the other hand we have
\begin{align*}
(e_{rs}\ot 1)E
&=\sum_{i\leq j} e_{rs}q_{ij} \ot p_{ij}\\
&=\sum_{i} e_{rs}e_{ii} \ot p_{ii} 
+\sum_{i<j} (e_{rs}e_{ii} +te_{rs}e_{ij})\ot p_{ij}\\
&=e_{rs}\ot p_{ss} + \sum_{s< j} (e_{rs}+ t e_{rj})\ot p_{sj}.
\end{align*} 
Similarly
\begin{align*}
E(e_{rs}\ot 1)
&=\sum_{i\leq j} q_{ij}e_{rs} \ot p_{ij}\\
&=\sum_{i} e_{ii} e_{rs}\ot p_{ii} 
+\sum_{i<j} (e_{ii}e _{rs} +te_{ij}e_{rs})\ot p_{ij}\\
&=e_{rs}\ot p_{rr} + \sum_{r<j} e_{rs} \ot p_{rj} + \sum_{i\leq r} t e_{is}\ot p_{ir}.
\end{align*} 
We see that $E(e_{rs}\ot 1)\notin C_u\ot P_u$ while $(e_{rs}\ot 1)E\notin C_u\ot P_u$.
\evoorb
For this example, we will get that the associated coproduct is a homomorphism. It has regular canonical maps $T_1$ and $T_3$. Again $T_2$ and $T_4$  will not be regular. We see that the right leg is all of $P_u$. The left leg of $E$ is spanned by the elements $e_{ii}$ and $e_{ii}+te_{ij}$ with $i<j$. Hence by all $e_{ij}$ with $i\leq j$. So $E$ is a full idempotent and $\Delta$ will be a full coproduct.
\ssnl
The next example is very similar.

\voorb\label{voorb:3.25}
We now define
Let $E=\sum_{i,j} q_{ij}\ot (p_i\ot p_j)$ where
\begin{equation*}
q_{jj} = e_{jj} \tussenen q_{ij}=e_{ii}+te_{ij} \text{ if } i\neq j.
\end{equation*}
We still will have $E(1\ot P)\subseteq C\ot P$ and $(1\ot P)E\subseteq C\ot P$ but now both $(C\ot 1)E$ and $E(C\ot 1)$ will no longer be subsets of $C\ot P$. Still $E$ is a full idempotent, it is even self-adjoint.
\ssnl
For the corresponding coproduct $\Delta$, we obtain that it is a $^*$-homomorphism and that it is full. The maps $T_1$ and $T_3$ are regular, but $T_2$ and $T_4$ are not.
\evoorb

It should be clear that, playing around with these ideas, one can construct several other examples of non-regular coproducts.
\snl
However, it is not completely clear how far one can get. One still has to find non-regular multiplier Hopf algebras with non-regular coproducts, non-regular separability idempotents and non-regular weak multiplier Hopf algebra with such a non-regular canonical idempotent.
\keepcomment{
\ssnl
Aankondigen en verwijzing maken naar het artikeltje \cite{VD-lumalg} en ook \cite{VD-nr}. Verwijzing in de conclusie sectie maken.
\ssnl
Het hangt er wat van af wat we daar uiteindelijk kunnen realiseren?
\ssnl
Niet zo duidelijk wat we hier bedoeld hebben.}{}

%
%

 \section{\hspace{-17pt}. Conclusions, more remarks and possible further research}\label{s:concl} 
 
The concept of a coalgebra is well-established. It is a vector space $A$ with a coproduct $\Delta$ and a counit $\varepsilon$. The coproduct is a linear map $\Delta$ from $A$ to $A\ot A$ satisfying coassociativity $(\Delta\ot\iota)\Delta=(\iota\ot\Delta)\Delta$ while the counit is a linear functional satisfying $(\iota\ot\varepsilon)\Delta(a)=a$ and $(\varepsilon\ot \iota)\Delta(a)=a$ for all $a\in A$. See e.g. \cite{A}, \cite{S} and \cite{R}.
\ssnl
Unfortunately, this notion turns out to be too restrictive for the study of coproducts on algebras without identity. This is seen and explained in Example \ref{voorb:KG} in the beginning of Section \ref{s:copr}. We need another notion when we want to study multiplier Hopf algebras or weak multiplier Hopf algebras and coactions for these objects.
\ssnl
We have explained the problem with the notion of coassociativity for such a more general concept and we have given a couple of possible workable solutions. The material is not completely new and many of the concepts and results have been treated in earlier work. See e.g.\ \cite{VD-mha}, \cite{VD-W0} and \cite{VD-W1}. Similarly, we have problems with defining coactions. We plan to treat these problems in a separate paper, see \cite{VD-co}.
\ssnl
The aim of this note is to provide some more detailed arguments that are not found in these papers. We also have clarified some minor problems that have been overlooked before.
\ssnl
\emph{Another problem} is that the dual space $A'$ can not be made into an associative algebra when a coproduct on $A$ maps into the multiplier algebra $M(A\ot A)$ and not into $A\ot A$. This problem is treated in Section \ref{s:dual} and several solutions are developed.
\ssnl
Finally, in the last section of this paper, we have given a few  examples. There are examples where some of the canonical maps are regular and others are not. We discussed the notions of coassociativity for such examples. For one specific example, we have also explained some of the results concerning the dual algebra construction as developed in Section \ref{s:dual}.
\ssnl
More research is welcome. One of the remaining problems is finding non-regular multiplier Hopf algebras and weak Hopf algebras. A first step towards such examples is finding non-regular separability idempotents. Some of the examples given in the last section of this paper provide certainly some ideas to find such examples.
\oldcomment{Hier mag nog wel wat aan toegevoegd worden - zie artikel5ra.tex}{}

%
%



\begin{thebibliography}{99}
\footnotesize \itemsep=0pt

\oldcomment{
\bibitem{??}  
\ssnl
 Op plaatsen waar nog een referentie moet gezocht worden - verschijnt als [1]
 \ssnl
}{}
\bibitem{A} E.\ Abe: {\it Hopf algebras}. \rm Cambridge University Press (1977).

\bibitem{Da} J.\  Dauns : {\it Multiplier rings and primitive ideals}. Trans.\ Amer.\ Math.\ Soc.\ 145 (1969), 125-158.

\bibitem{De-VD} L.\ Delvaux \& A.\ Van Daele: {\it Algebraic quantum hypergroups}.  Adv.\ Math.\ 226 (2011), 1134-1167. 
See also Arxiv math.RA/0606466.

\bibitem{La-VD1} M.B.\ Landstad \& A.\ Van Daele: {\it Duality of algebras with an antipode and integrals}. Preprint  NTNU (Trondheim, Norway) and KU Leuven (Leuven, Belgium). In preparation. 

\bibitem{P} G.K.\ Pedersen: {\it C*-algebras and their automorphism groups}. Academic Press, New York (1979).

\bibitem{R} D.\ Radford: {\it Hopf algebras}. Series on Knots and Everything Vol. 49, Word Scientific, Singapore (2012).

\bibitem{S} M.\ Sweedler: {\it Hopf algebras}. Benjamin, New-York (1969).

\bibitem{T-VD} T.\ Timmermann \& A.\ Van Daele: {\it Regular multiplier Hopf algebroids. Basic theory and examples}.
Commun. in Alg. 46 (2017), 1926-1958. 

\bibitem{T-VD2} T.\ Timmermann, A.\ Van Daele \& S.Wang: {\it Pairing and duality of algebraic quantum groupoids.}
Int.\ J.\ Math.\ Vol.\ 33 (2022), 1–50.

\bibitem{VD-mha} A.\ Van Daele: {\it Multiplier Hopf algebras}. Trans. Am. Math. Soc.  342(2) (1994), 917-932.

 \bibitem{VD-alg} A.\ Van Daele: {\it An algebraic framework for group duality}. Adv. in Math. 140 (1998), 323-366.

\bibitem{VD-co} A.\ Van Daele: {\it Reflections on coactions for non-unital algebras}. Preprint KU Leuven (Leuven, Belgium). In preparation.

\bibitem{VD-lumalg} A.\ Van Daele \& J.\ Vercruysse: {\it Multiplier algebras and local units}. Preprint KU Leuven (Leuven, Belgium) and VUB (Brussels, Belgium). In preparation. \oldcomment{Titel ok?}{} 

\bibitem{VD-infmalg} A.\ Van Daele  \& J.\ Vercruysse: {\it Infinite matrix algebras}. Preprint University KU  Leuven (Leuven, Belgium) and VUB (Brussels, Belgium).  In preparation. \oldcomment{Titel ok?}{} 

\bibitem{VD-W0} A.\ Van Daele \& S.\ Wang: {\it Weak multiplier Hopf algebras. Preliminaries, motivation and basic examples}. Operator Algebras and Quantum Groups. Banach Center Publications 98 (2012), 367-415. 

\bibitem{VD-W1} A.\ Van Daele \& S.\ Wang: {\it Weak multiplier Hopf algebras I. The main theory}. Journal f\"ur die reine und angewandte Mathematik (Crelles Journal) 705 (2015), 155-209, ISSN (Online) 1435-5345, ISSN (Print) 0075-4102, DOI: 10.1515/crelle-2013-0053, July 2013. 

\end{thebibliography}
\end{document}